\documentclass[12pt,amstex]{article}
\usepackage[usenames]{color}

\usepackage{amscd}
\usepackage{amsmath}
\usepackage{amssymb}
\usepackage{graphicx}

\newtheorem {theo} {\bf Theorem} [section]
\newtheorem {prop} [theo] {\bf Proposition}

\newtheorem {lem} [theo] {\bf Lemma}
\newtheorem {defn} [theo] {\bf Definition}

\newtheorem {rem} [theo] {\bf Remark}
\newcommand{\QED}{\hfill \CaixaPreta \vspace{6mm}}
\def\CaixaPreta{\vrule Depth0pt height6pt width6pt}

\newcommand{\qed}{\nopagebreak\hfill{\vrule width6pt height6pt depth0pt}}
\newcommand{\be}{\begin{eqnarray}}
\newcommand{\ee}{\end{eqnarray}}
\newcommand{\benn}{\begin{eqnarray*}}
\newcommand{\eenn}{\end{eqnarray*}}
\newcommand{\bse}{\begin{equation}}
\newcommand{\ese}{\end{equation}}
\newcommand{\bsenn}{\begin{displaymath}}
\newcommand{\esenn}{\end{displaymath}}
\newcommand{\logand}{\;\;{\rm and }\;\;}

\newcommand{\logor}{\;\;{\rm or }\;\;}
\newcommand{\logif}{\;\;{\rm if }\;\;}

\newcommand{\such}{\;\;{\rm such\; that }\;\;}

\newcommand{\R}{\mathbb{R}}

\newcommand{\ke}{\mathrm{Ker}\,}
\newcommand{\im}{\mathrm{Im}\,}
\newcommand{\rank}{\mathrm{Rank}\,}
\newcommand{\nullity}{\mathrm{Nullity}\,}
\newcommand{\Ln}{\mathrm{Ln}\,}
\newcommand{\Exp}{\mathrm{Exp}\,}
\usepackage{graphicx}
\usepackage{url}
\usepackage{amsmath}

\begin{document}
\title{Chemical Reaction Networks in a Laplacian Framework}
\author{J. J. P. Veerman\thanks{Weizmann Institute of Science, Rehovot, Israel}
	\thanks{Fariborz Maseeh Dept. of Math. and Stat., Portland State Univ., Portland, OR, USA; e-mail: veerman@pdx.edu}, T. Whalen-Wagner\thanks{Fariborz Maseeh Dept. of Math. and Stat.,
Portland State Univ., Portland, OR, USA; e-mail: tessaww@hushmail.com}, Ewan Kummel\thanks{Fariborz Maseeh Dept. of Math. and Stat.,
Portland State Univ., Portland, OR, USA; e-mail:ewan@pdx.edu}\\
}\maketitle

\begin{abstract}
The study of the dynamics of chemical reactions, and in particular phenomena such as
oscillating reactions, has led to the recognition that many dynamical properties of a chemical reaction
can be predicted from graph theoretical properties of a certain directed graph, called
a Chemical Reaction Network (CRN). In this graph, the edges represent the reactions and the vertices
the reacting combinations of chemical substances.

In contrast with the classical treatment, in this work, we heavily rely on a recently developed theory
of directed graph Laplacians to simplify the traditional treatment of the so-called deficiency zero
systems of CRN theory. We show that much of the dynamics of these polynomial systems of
differential equations can be understood by analyzing the directed graph Laplacian associated
with the system. Beside the more concise mathematical treatment, this leads to considerably stronger
results. In particular, (i) we show that our Laplacian deficiency zero theorem is markedly stronger than
the traditional one and (ii) we derive simple equations for the locus of the equilibria in all
(Laplacian) deficiency zero cases.

This paper is written in a way to make the material easily accessible to a mathematical audience.
In particular, no knowledge of chemistry or physics is assumed.
\end{abstract}


\vskip .0in

 %\normalsize        %\mysetfontsize5

%\vskip 0.4in
\begin{centering}\section{Introduction}
 \label{chap:intro}\end{centering}
\setcounter{figure}{0} \setcounter{equation}{0}

CRN's form a compelling area of study with many connections to other areas of mathematics. For example,
recently, a wonderful introduction appeared highlighting the connection with algebraic geometry
\cite{dickenstein_2021}. In this paper, we review the basic theory of Chemical Reaction Networks (CRN's)
employing the recently developed formalism of directed graph Laplacians \cite{chebotarev2002, caugh,
mirzaev_2013, kummel, lyons}. In the literature since the 1970's
\cite{horn_1972, horn_jackson_1972, feinberg_1972}, this analysis has been based (for a variety of reasons)
on the understanding of a different linear operator that, however, contains less information.
The change to a Laplacian formulation allows us to give
a more concise derivation of all the classical results of the zero deficiency theory with much less
effort. Equally important is the fact that the Laplacian formulation gives stronger results,
as we explain below. With this paper, we wish to make the material accessible to a
mathematical audience. Thus we restrict our vocabulary to terms current in mathematics
or at least mathematical graph theory.

Here is an overview of what we aim to achieve in this work. A chemical reaction network (or CRT) consists
of a (often very large) collection of  first order polynomial differential equations. First we formulate
the \emph{Laplacian} version of the zero deficiency condition (Definition \ref{def:Ldeficiency}) which
essentially eliminates `unobserved' chemical reactions. We assume this condition for the rest of the paper.
We then prove the zero \emph{Laplacian} deficiency theorem (Theorem \ref{thm:0deficiency}), which says
that there exists a strictly positive equilibrium if and only if the associated directed graph is
componentwise strongly connected  (or CSC, see Definition \ref{def:connected}). Subsequently, we will
prove that for every choice of certain constants of the motion (Definition \ref{def:affinespace}), there
is exactly one equilibrium (Theorem \ref{thm:EandU}) and furthermore that this equilibrium is locally
asymptotically stable (Theorem \ref{thm:convergence}).

As mentioned, the Laplacian framework allows us not only to give more concise proofs but also leads to
stronger results. We now describe the new aspects of this work.
The Laplacian zero deficiency theorem is strictly stronger than its classical counterpart (Proposition
\ref{prop:USimpliesTHEM2}) and in Section \ref{chap:Comparison1} we give a significant example of that.
Using the Laplacian theory, we can in fact show that the existence of a positive equilibrium
in any zero deficiency system is equivalent
to the existence of an orbit in a compact subset of the open positive orthant (Theorem
\ref{thm:0deficiencybdd}). For general zero deficiency systems, we derive simple equations that
determine the locus of any equilibrium (Theorem \ref{thm:0deficiencyeqn}).
We give examples of this
in Section \ref{chap:examples}. Finally, in some cases, the Laplacian method detects more constants of
the motion than the traditional one (Proposition \ref{prop:USimpliesTHEM1} and the remarks that follow it).

The first sporadic accounts of oscillating chemical reactions were published in the 19th century.
At the time, they received very little attention, in part because known examples were difficult to
reproduce and in part because of a belief among scientists that such behavior was impossible. When
Bray \cite{bray_1921} published the first detailed description of such a reaction in the 1920's, the
consensus among his peers was that the behavior must be the result of experimental error.
Indeed, 30 years later Belousov spent 8 years trying to publish a description of his famous reaction.
His observations were eventually published in a non-peer reviewed journal (see \cite{winfree_1984}
for details). Belousov's publication
allowed other researchers to replicate his example, produce others, and eventually derive conditions
needed for such reactions \cite{zhabot_1991}.

While the study of chemical reaction networks is at least old as the introduction of detailed balance
for chemical reactions \cite{wegsch_1901}, the mathematical theory of chemical reaction networks began
in earnest in the 1960's with the work of Aris \cite{aris_1965} and
achieved prominence in the 1970's with the work of Horn, Jackson, Feinberg, and others
\cite{feinberg_1972, horn_jackson_1972, horn_1972, feinberg_horn_1974}, see also \cite{feinberg_2019}.
Since reaction rates are difficult to
measure experimentally, this theory was in part motivated by the need to understand exotic
behaviors of chemical reactions in a way that does not require knowing precise reaction rates.
Examples of such behavior are oscillations and bi-stability.
The landmark 1987 Feinberg paper \cite{feinberg_1987} combines much of this early work into two theorems:
the deficiency zero theorem which we discuss below, and an extension called the deficiency one
theorem. In each case, the deficiency (Definition \ref{def:Ldeficiency}) of a reaction network
is used to characterize the equilibria.

Probably the first papers in which a chemical reaction diagram, consisting of chemical compounds
connected by arrows, was explicitly treated as a directed graph was \cite{sontag-2001, sontag-2002}. It turns
out that many notions in the theory of chemical reaction networks have direct parallels in the
language of directed networks. Indeed, strictly from the perspective of dynamics, there is
independent mathematical interest in the notion that the behavior of the highly non-linear system
can be at least partly understood by the analysis of directed graphs. There are many recent papers
\cite{schaft_2013, schaft_rao_jayaw_2015, kim_2018} that make use of this interplay to derive new
results. Other work \cite{gunawardena_2003, craciun_dickenstein_shiu_sturmfels_2009}
specifically exploited the structure of directed graph Laplacians. However, they did
so without the benefit of a clear, standardized theory describing such Laplacians.

This is an area of active on-going research. One focus of research is the global attractor conjecture, which
asserts that if the associated directed graph is componentwise strongly connected in the zero deficiency case,
 then \emph{every}
initial condition in the open orthant converges to an equilibrium (see the remark after Theorem \ref{thm:convergence}).
This conjecture appears as early as \cite{horn_jackson_1972}, where it was mistakenly believed
to be proved, and has been shown in certain cases \cite{anderson_2011}. The analysis of higher deficiency reaction
networks is another active area of study \cite{joshi_2017, pucci_2018,kaufman_soule_2019}. Some
results can be extended fairly easily to the deficiency one case as discussed in  \cite{feinberg_2019},
but in general the behavior of higher deficiency networks is not well understood.
Another important open question is that of
``persistence". In general, a persistent reaction network is one in which all chemical concentrations
have a positive lower bound for all positive time \cite{brunner_2018, feinberg_2019}.
One famous and, so far, unproved conjecture is that this holds for every network whose associated
directed graph is componentwise strongly connected (independent of the deficiency) \cite{shiu_2010}.
See also the comment after Theorem \ref{thm:0deficiencybdd}.

The original emphasis in reaction networks was biased towards controlling chemical reactions and
therefore trying to ensure that exotic behavior does not arise. However, in recent years the promise
of applying CRN theory to complex biological systems has shifted that interest towards seeking out
and analyzing more complicated behaviors \cite{gunawardena_2003}. For instance, for the high deficiency
case, it possible that reactions take place even though the associated linear system of reaction
equations is at an equilibrium.  For a collection of examples, we refer the reader to
\cite{feinberg_2019}.

The outline of this paper is as follows. We first (Section \ref{chap:prelim}) discuss some well-known
preliminary results that we will need later, as well as some notation. In Section \ref{chap:Laplacians},
we summarize the modern theory of directed graph Laplacians and its conclusions. Section \ref{chap:chem}
describes the mathematical definition of chemical reaction networks. Section \ref{chap:0defic} states and
proves the zero deficiency theorem. In Sections \ref{chap:uniqueness} and \ref{chap:Lyapunov}, we prove
that in zero deficient systems satisfying a certain connectedness property, each invariant subspace
has a unique asymptotically stable equilibrium. In Section \ref{chap:examples}, we give
a few examples of reaction networks designed to illustrate the theory. In Section \ref{chap:Comparison1},
we compare our results and their classical counterparts and show that our results in some cases
improve classical results.

\noindent
{\bf Acknowledgement:} We are grateful to Patrick de Leenheer and Arjan van der Schaft for helpful
conversations. We also wish to thank the referee for the valuable comments which improved
the paper substantially.

\begin{centering}\section{Preliminaries}
 \label{chap:prelim}\end{centering}
\setcounter{figure}{0} \setcounter{equation}{0}

In this section, we summarize some well-known results that we will need to use later and present
some notation. The first two lemmas are standard results of linear algebra. Let $A:\R^n\rightarrow \R^m$
and $B:\R^e\rightarrow \R^n$ be linear maps and $V$ and $W$ (linear) subspaces of $\R^n$.

\begin{lem} For linear subspaces $V$ and $W$: $V^\bot\cap W^\bot=(V+W)^\bot$.
\label{lem:linalg1}
\end{lem}

\begin{lem} For any matrix $A$ we have: $\ke  A=(\im  A^T)^\bot$, where the orthogonal complement
is in the domain of $A$.
\label{lem:linalg2}
\end{lem}

\vskip-0.1in\noindent
It follows that $\ke  A$ and $\im  A^T $ span the domain of $A$ and so the sum of their dimensions equals $n$.
$\dim \ke A$ is referred to as the nullity of $A$ and $\dim \im A^T$ is equal to the rank of $A$.

Putting the previous lemmas together, we immediately see the following.

\begin{prop} For any two matrices $A$ and $B$: $\left[\ke  A\cap \im  B \right]^\bot =
\im  A^T + \ke  B^T$.
\label{prop:linalg4}
\end{prop}

\begin{prop} For any two matrices $A$ and $B$:
$\dim \left[\ke  B \cap \im A \right]= \dim \ke BA - \dim \ke A$.
\label{prop:linalg5}
\end{prop}

\vskip0.0in\noindent
{\bf Proof.} To prove the equality, it is sufficient to show that the linear map $x\rightarrow Ax$
induces a bijection
\bsenn
\psi:\ke BA/\ke A\rightarrow \ke A \cap \im B \,.
\esenn
Indeed, $\psi$ is well-defined and injective, because for $x$ and $y$ in $\ke BA$:
\bsenn
Ax=Ay \quad \Longleftrightarrow \quad A(y-x)=0 \quad \Longleftrightarrow \quad y-x \in \ke A \,.
\esenn
Clearly $\psi$ is surjective, because for any $z\in \ke B \cap \im A$, there is an $x$ such that $z=Ax$.
\QED

\vskip-0.0in
We will also need a few simple calculus lemmas.

\begin{lem} For any $a>0$ and $b>0$, we have: $a(\ln a -\ln b)\geq a-b$. Equality iff $a=b$.
\label{lem:calc1}
\end{lem}

\vskip0.0in\noindent
{\bf Proof.} The tangent line to $\ln x$ at $x=1$ is above the graph of that function for all $x\neq 1$,
and so $x-1\geq \ln x$. Substituting $x=b/a$ yields the result. \QED

\vskip-0.2in\noindent
\begin{lem} For any $a>0$ and $b>0$, we have: $(a-b)(\ln a -\ln b)\geq 0$.
Equality iff $a=b$.
\label{lem:calc2}
\end{lem}

\vskip0.0in\noindent
{\bf Proof.} Lemma \ref{lem:calc1} implies $-b(\ln a -\ln b)\geq b-a$. Adding that inequality
to the one in Lemma \ref{lem:calc1} proves the result. \QED

\vskip-0.2in\noindent
\begin{lem} For any $x>0$ and $z>0$, there are $\mu_\pm \in \R$ such that
$\forall \mu \not\in [\mu_-,\mu_+]\;:\; xe^{\mu}-z\mu>x$.
\label{lem:calc3}
\end{lem}

\vskip0.0in\noindent
{\bf Proof.} Taking the derivative of $f(\mu):= xe^{\mu}-z\mu$ shows that this function has a global
minimum at $\mu^*=\ln z - \ln x$. Applying Lemma \ref{lem:calc1} to $f(\mu^*)=z-z(\ln z-\ln x)$ shows
that $f(\mu^*)\leq z$.
Finally, $\lim_{\mu\to \pm \infty} f(\mu)=+\infty$. \QED

\vskip-0.2in%\noindent
Finally, we need a result from the theory of dynamical systems.

\vskip-0.1in\noindent
\begin{defn} A function $V:{\cal O}\subseteq \R^n \rightarrow \R$ where ${\cal O}$ is open, is called
a Lyapunov function for the system $\dot x=f(x)$ in $\R^n$ if it is a continuously differentiable
and satisfies that along a trajectory $\dot V(x(t)):=(\nabla V(x(t)),\dot x(t)) \leq 0$, where $(,)$ is
the standard inner product, and $\nabla$ the gradient, both in $R^n$.
\label{def:lyapunov}
\end{defn}

\vskip-0.1in\noindent
\begin{defn} The $\omega$-limit set of $x$ is the set of points $y$ for which there is a sequence
$t_n\rightarrow \infty$ so that $\lim_{n\rightarrow \infty}\,x(t_n) = y$
\label{def:omegalimit}
\end{defn}

\begin{theo}{\cite{stern}} Let $V:{\cal O}\rightarrow \R$ be a Lyapunov function for the system $\dot x= f(x)$.
The intersection of the $\omega$-limit set of a point $x$
and the set ${\cal O}$ is contained in the set where $(\nabla V(x(t)),\dot x(t))=0$.
\label{thm:lyap}
\end{theo}

\vskip-0.1in\noindent
From now on, we will use the abbreviation $\dot V$ instead of the cumbersome $(\nabla V(x(t)),\dot x(t))$.

Finally, we mention some notation that will be used throughout this paper.
Given vectors $x$ and $y$ in $\R^n$, we define $x\odot y$ as the vector whose components are $x_iy_i$.
(This is also called the Hadamard product.)
We write $x/y$ for the vector whose components are $\frac{x_i}{y_i}$ ($y_i\neq 0$ for all $i$).
The componentwise logarithm of $x$  ($x_i > 0$ for all $i$) is denoted by $\Ln x$, while
the componentwise exponential of $x$ will be written as $\Exp x$.  We write $x>0$ when
$x_i > 0$ for all $i$. Given a system of differential equations $\dot x=f(x)$, we will use the word \emph{equilibrium} for a point $x$ such that $f(x)=0$.

%\vskip 0.4in
\begin{centering}\section{Laplacians}
 \label{chap:Laplacians}\end{centering}
\setcounter{figure}{0} \setcounter{equation}{0}

\begin{figure}[!ht]
\center
\includegraphics[width=4.0in,height=1.6in]{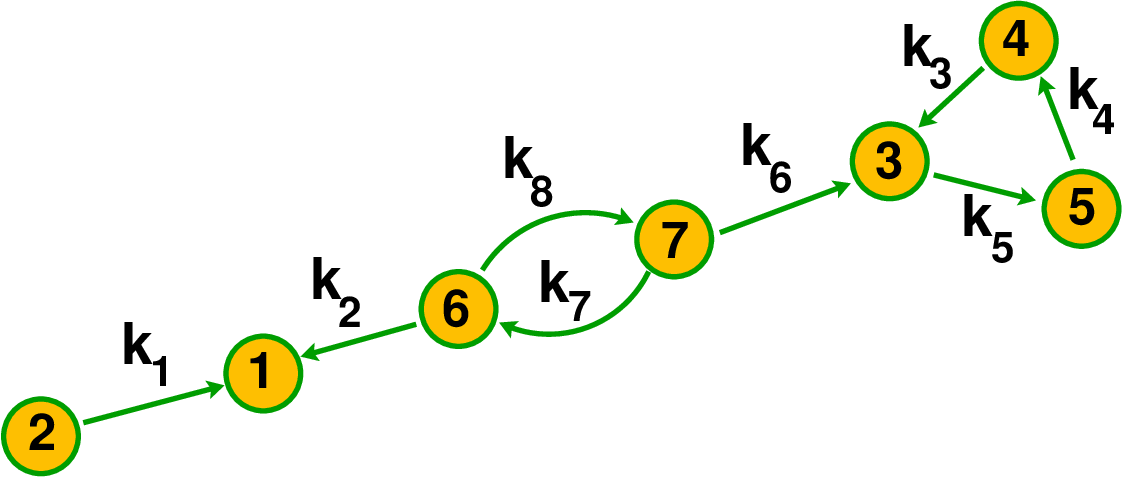}
\caption{\emph{An example of a directed network $G$. See also Section \ref{chap:examples}.}}
\label{fig:example-reversed}
\end{figure}

\begin{figure}[!ht]
\center
\includegraphics[width=4.0in,height=1.6in]{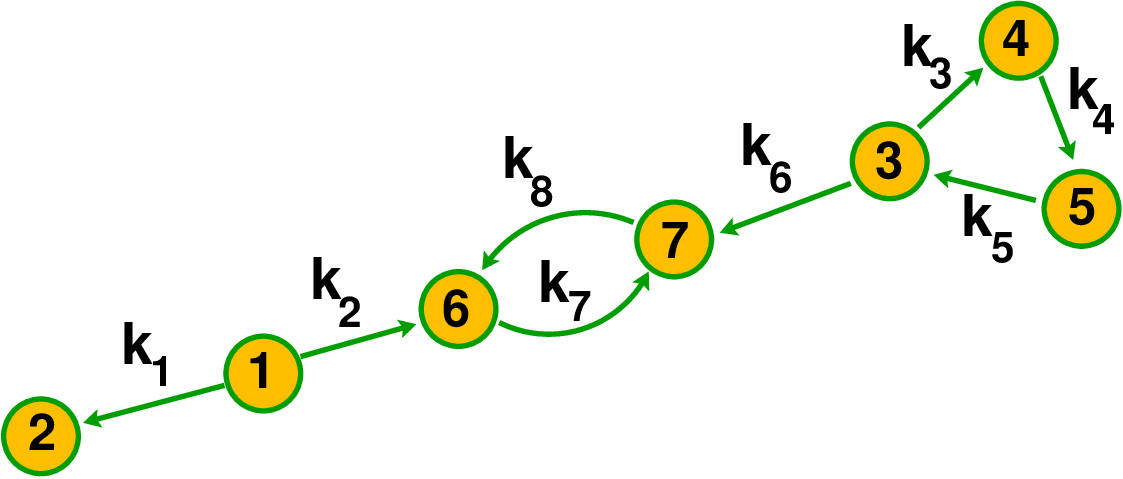}
\caption{\emph{The directed network $G^\Lsh$ obtained from Figure \ref{fig:example-reversed}
by reversing the orientation of the edges.}}
\label{fig:example}
\end{figure}

Two things are important to bear in mind when working with \emph{directed} graphs to model certain phenomena.
First of all, directed graphs are used to model interactions that are not symmetric, i.e. the influence
of $x$ on $y$ may not be the same as the influence of $y$ on $x$. As a result, the Laplacian
is usually not symmetric and its eigenvalues are not necessarily real. Another complicating factor
is that different authors may choose opposite orientations of the edges. Below, we will use $G^\Lsh$
for the graph obtained from a graph $G$ by reversing the orientation of all edges (compare Figures
\ref{fig:example-reversed} and \ref{fig:example}).

The conventions outlined in this section are taken from \cite{caugh, kummel, lyons}.

We now give a few of the basic facts of Laplacian dynamics in a directed, loopless\footnote{A loop is an edge
that starts and ends at the same vertex} graph $G$. We assume that $G$ has $v$ vertices and $e$ directed edges.
The $v\times e$ matrix $B$ is the \emph{begin matrix} \cite{DG4} such that $B_{ij}=1$ if vertex $i$ starts
edge $j$ and 0 otherwise. Similarly, the $v \times e$ \emph{end matrix} \cite{DG4} $E$ is defined by
$E_{ij}=1$ if vertex $i$ ends edge $j$ and 0 otherwise. We use these matrices to define the
boundary operator (or \emph{incidence matrix} in graph theory texts) $\partial:=E-B$. As an example, we
exhibit the boundary operator associated the graph in Figure \ref{fig:example-reversed}:
\bse
\partial = E-B =
\begin{pmatrix} 1&1&0&0&0&0&0&0\\-1&0&0&0&0&0&0&0\\0&0&1&0&-1&1&0&0\\0&0&-1&1&0&0&0&0\\
0&0&0&-1&1&0&0&0\\0&-1&0&0&0&0&1&-1\\0&0&0&0&0&-1&-1&1 \end{pmatrix}
\label{eq:boundary-op}
\ese
The \emph{weight matrix} $W$ is diagonal with (strictly) positive weights on the diagonal.
The weights are equal to 1 in the \emph{unweighted} case.

\begin{defn}{\cite{DG4}} The \emph{undirected} weighted Laplacian $L$ (also called the Kirhhoff matrix)
is given by:
\bsenn
L\equiv \partial\,W \partial^T=(E-B)W(E-B)^T.
\esenn
It is the sum of the \emph{in-degree} Laplacian $L_{\mathrm{in}}$ and the \emph{out-degree} Laplacian
$L_{\mathrm{out}}$.
\bsenn
L_{\mathrm{in}}= EW(E-B)^T \quad \logand \quad L_{\mathrm{out}}= -BW(E-B)^T.
\esenn
Note that the out-degree Laplacian of $G$ is the same as the in-degree Laplacian of $G^\Lsh$.
\label{def:laplacians}
\end{defn}

As an example, we give the unweighted in- and out-degree Laplacian of the graph $G$ in Figure
\ref{fig:example-reversed}. The $i$th row of $L_{in}$ gives the vertices with edges coming \emph{to}
the $i$th vertex, while the the $i$th row of $L_{out}$ gives the vertices with edges coming \emph{from}
the $i$th vertex.
\be
L_{\textrm{in }} &=& \left( \begin{array}{ccccccc}
{{\bf 2}}&{{\bf -1}}&0&0&0&{{\bf -1}}&0\\
0&0&0&0&0&0&0\\
0&0&{{\bf 2}}&{{\bf -1}}&0&0&{{\bf -1}}\\
0&0&0&{{\bf 1}}&{{\bf -1}}&0&0\\
0&0&{{\bf -1}}&0&{{\bf 1}}&0&0\\
0&0&0&0&0&{{\bf 1}}&{{\bf -1}}\\
0&0&0&0&0&{{\bf -1}}&{{\bf 1}}
\end {array} \right) \\[0.4cm]
L_{\textrm{out}} &=& \left( \begin{array}{ccccccc}
0&0&0&0&0&0&0\\
{{\bf -1}}&{{\bf 1}}&0&0&0&0&0\\
 0&0&{{\bf 1}}&0&{{\bf -1}}&0&0\\
 0&0&{{\bf -1}}&{{\bf 1}}&0&0&0\\
 0&0&0&{{\bf -1}}&{{\bf 1}}&0&0\\
 {{\bf -1}}&0&0&0&0&{{\bf 2}}&{{\bf -1}}\\
0&0&{{\bf -1}}&0&0&{{\bf -1}}&{{\bf 2}}
\end{array} \right)
\label{eq:indegree-outdegree}
\ee

\noindent
{\bf Remark:} More generally, a Laplacian is a square matrix with non-negative diagonal and non-positive
off-diagonal elements whose row-sums all give zero. It is easy to see that any such matrix can be written
in the form stipulated by Definition \ref{def:laplacians}.

\begin{defn}{\cite{caugh, kummel, lyons}} Given a directed graph $G$.\\
1) A vertex $j$ is in the reachable set from the vertex $i$ if $j=i$ or there is a directed path from
$i$ to $j$, $i\rightsquigarrow j$. The reachable set from $i$ is called $R(i)$. \\
2) A \emph{reach} $R$ is a maximal reachable set (i.e. one that is not \emph{properly} contained in any other reachable set $R(j)$).\\
3) A \emph{cabal} $C\subseteq R$ is the maximal (largest) set of vertices from which all of a reach $R$
is reachable.\\
4) The exclusive part $H$ of a reach $R$ is the set of vertices contained in $R$ and
in no other reach.\\
5) The common part $C = R - H$ is the set of vertices that $R$ has in common with some other reach.
\label{def:reach}
\end{defn}

To illustrate these notions, consider the network of Figure \ref{fig:example-reversed}.
There are two reaches: $R_1=\{2,1\}$ and $R_2=\{1,6,7,3,4,5\}$.
Their cabals are given by $B_1=\{2\}$ (in $R_1$) and $B_2=\{6,7\}$ (in $R_2$).
(We note that a cabal consisting of a single vertex as is the case in $R_1$ is often called a
\emph{leader}.)
In $R_1$, only $H_1=\{1\}$ is not shared by $R_2$, and so $C_1=\{2\}$. Similarly,
$H_2=\{6,7,3,4,5\}$, while $C_2=\{1\}$.

\begin{defn}{\cite{caugh, kummel, lyons}}
A \emph{co-reach} is a reach in $G^\Lsh$ and a \emph{co-cabal} is a cabal in $G^\Lsh$.
\label{def:coreach}
\end{defn}

Thus the \emph{co-reaches} and \emph{co-cabals} of $G$ in Figure \ref{fig:example-reversed}
can be found as the reaches and cabals of $G^\Lsh$ in Figure \ref{fig:example}.
They are given by $R_1^\Lsh=\{2,1,6,7\}$ with cabal $B_1^\Lsh=\{1\}$ and $R_2^\Lsh=\{6,7,3,4,5\}$
with cabal $B_2^\Lsh=\{3,4,5\}$.

In the following, the right kernel of a matrix $A$ denotes the set of vectors $x$ so that $Ax=0$,
while the left kernel is the set of (row) vectors $y$ such that $yA=0$.

\begin{theo}{\cite{caugh, kummel, lyons}} Let $G$ be a digraph with reaches $R_1$,..., $R_k$. The eigenvalue 0
of $L_{\mathrm{in}}$ has geometric and algebraic multiplicity $k$. All other eigenvalues have negative
real part.
\label{thm:multiplicity}
\end{theo}

\vskip -0.2in \noindent
\begin{theo}{\cite{caugh, kummel, lyons}} Let $G$ be a digraph with reaches $R_1$,..., $R_k$. The
\emph{column} vectors $\{\gamma_1,\cdots, \gamma_k\}$ form a basis for the \emph{right} kernel of a
Laplacian $L_{\mathrm{in}}$, where:
\bsenn
\left\{\begin{matrix}
\gamma_{m,j}=1 & \logif & j\in H_m &\textrm{(exclusive)}\\
\gamma_{m,j} \in (0,1) & \logif & j\in C_m &\textrm{(common)}\\
\gamma_{m,j}=0 & \logif & j\not\in R_m &\textrm{(not in reach)}\\
\sum_{m=1}^k\,\gamma_{m,j} ={\bf 1} &   &  &
\end{matrix} \right.
\esenn
\label{thm:rightkernel}
\end{theo}

\vskip -0.4in \noindent
\begin{theo}{\cite{kummel, lyons}} Let $G$ be a digraph with $k\geq 1$ reaches. The \emph{row} vectors
$\{\bar\gamma_1,\cdots, \bar\gamma_k\}$ form a basis for the \emph{left} kernel of a Laplacian
$L_{\mathrm{in}}$, where:
\bsenn
\left\{\begin{matrix}
\bar\gamma_{m,j}>0 & \logif & j\in B_m & \textrm{(cabal)}\\
\bar\gamma_{m,j}=0 & \logif & j\not\in B_m & \textrm{(not in cabal)}\\
\sum_{j=1}^k\,\bar\gamma_{m,j}= 1 &   & &
\end{matrix} \right.
\esenn
\label{thm:leftkernel}
\end{theo}

\vskip -0.3in \noindent
Later on, it will be of considerable importance that the vectors $\{\bar\gamma_1,\cdots, \bar\gamma_k\}$
have disjoint support (as opposed to the vectors $\{\gamma_1,\cdots, \gamma_k\}$) and so form an
orthogonal basis of the left kernel. We note in passing that these null vectors form a basis for the
space of stationary distributions in Markov chains \cite{kummel, lyons} and are related to the maximal
spanning forests \cite{chebotarev2002}.

\vskip -0.1in \noindent
\begin{defn} For directed graphs, we distinguish \emph{weakly connected components} -- a maximal set
of vertices for which there is an \emph{undirected} path between every pair of vertices -- from
\emph{strongly connected components} (SC's) -- a maximal group of vertices for which there
is a \emph{directed} path between every pair of vertices.
\label{def:connected}
\end{defn}

\vskip -0.2in \noindent
\begin{defn} A graph $G$ is componentwise strongly connected (abbreviated to CSC) if every
weak component is a strong component.
\label{def:CSC}
\end{defn}

\vskip -0.1in \noindent
{\bf Remark:} One easily sees that the following statements are equivalent: \\
$G$ is CSC,\\
$G^\Lsh$ is CSC, \\
every reach of $G$ (or $G^\Lsh$) is strongly connected, and\\
every reach is a cabal.

\vskip-0.0in\noindent
\begin{lem} For any $G$, $\rank \partial$ (or $\dim \im \partial$) equals the number of vertices minus
the number of weak components. Furthermore, $\nullity \partial^T$ (or $\dim \ke \partial^T$) equals the
number weak components.
\label{lem:dimdelta}
\end{lem}

\vskip0.0in\noindent
{\bf Proof.} It is sufficient to prove this if $G$ consists of one weak component.
Suppose $G$ is a weak component with $v$ vertices and $k$ reaches $\{R_i\}_1^k$.
Note that since every edge has one endpoint and one begin point, $\partial^T {\bf 1}=0$.
Lemma \ref{lem:linalg2} and the remark following it now imply $\nullity \partial^T + \rank \partial = v$. Thus
$\rank \partial\leq v-1$.

Choose an arbitrary vertex $b$ as basepoint and let $r\ne b$ be any other vertex.
By definition \ref{def:connected}, there is an undirected path $\gamma$ from $b$ to $r$.
Consider the directed edges $e_i$ of $\gamma$. When $\gamma$ traverses $e_i$ in the positive
direction, multiply $e_i$ by $w_i=+1$, and in the other case by $w_i=-1$. The image under $\partial$ of
$\sum w_ie_i$  is $r-b$. This shows that $\rank \partial\geq v-1$. The second statement follows
from the first, because the remark after Lemma \ref{lem:linalg2} says that the sum of the two dimensions
must be $v$.
\QED

\vskip-0.1in\noindent
{\bf Remark:} This lemma is standard fare in algebraic graph theory \cite{bollobas}
and, in fact, algebraic topology. In a nutshell,
it is how the zeroth homology and cohomology groups are computed. To illustrate the procedure
in the second paragraph of the proof, we turn to Figure \ref{fig:example-reversed}.
In that figure, denote the vertices marked with $i$ by $v_i$ and the edges marked with $k_j$ by $e_j$.
Choose, for example, basepoint $b=v_1$ and endpoint $r=v_7$. A path $\gamma$ from $b$ to $r$ is given by
$-e_2+e_8$. Apply the boundary operator (read off from Figure \ref{fig:example-reversed} or using
\eqref{eq:boundary-op}) to get $\partial\gamma=-(v_1-v_6)+(v_7-v_6)=v_7-v_1$.

%\vskip 0.4in
\begin{centering}\section{Chemical Reaction Networks with Mass Action}
 \label{chap:chem}\end{centering}
\setcounter{figure}{0} \setcounter{equation}{0}

The three basic ingredients of a CRN are:
\bsenn
\left\{\begin{matrix} c \textrm{ ``concentrations of molecules or similar chemical substances", each
denoted by } x_i \,;\\
                      v \textrm{ vertices or ``concentrations of reacting mixtures", each denoted by }v_i\,;\\
                      e \textrm{ directed edges or ``reaction rates", each denoted by }e_i\,.\end{matrix} \right.
\esenn
We then associate a linear vector space to each of these ingredients as follows.
The column vectors $(x_1,\cdots, x_c)^T$ form the space $\R^{c}$. In the same way, vectors in
the space $\R^v$ and $\R^e$ have components $v_i$ and $e_i$, respectively.
The begin and end matrices $B$ and $E$ defined in Section
\ref{chap:Laplacians} correspond to linear transformations $B$ and $E$ from $\R^e$ to $\R^v$, whereas
their transpose acts in the opposite direction.

The spaces $\R^v$ and $\R^e$ are used to compute rates of change of concentrations, not the concentrations
themselves. As an example, look at the simple system consisting of the reaction $2H_2+O_2\rightarrow 2H_2O$.
The $x_i$ are the concentrations of, respectively, $H_2$, $O_2$, and $H_2O$. There are 2 vertices,
$v_1$ denotes the concentration of the combination $2H_2+O_2$ and $v_2$ that of $2H_2O$ and one edge
(or reaction) $v_1{\buildrel e_1 \over \rightarrow }\, v_2$. While the concentration
of $2H_2+O_2$ is an ambiguous concept, the rate of change of that same quantity due to the reaction, is not.

Next, we describe the relationship between the reacting mixtures and the molecules. (Note that we are
dropping the quotation marks.) The count of $i$-molecules in the $j$th vertex -- or reacting mixture --
equals $S_{ij}$. Put more simply, the $j$th column of $S$ gives the composition of molecules in the $jth$
vertex.
Labeling both from left to right, the matrix $S$ for the reaction given above, is:
\bsenn
S=\begin{pmatrix} 2&0\\1&0\\0&2\end{pmatrix} \,.
\esenn
This defines a linear transformation $S:\R^v\rightarrow \R^c$ whose matrix has entries that are
\emph{non-negative integers}. In a system with many simultaneous reactions, the rate of change in $x_i$
(indicated by $\dot x_i$) equals the sum of the rates of change of those mixtures in which that molecule
occurs. Thus
\bse
\dot x=S\dot v \quad \logor \quad \dot x_i=\sum_j\,S_{ij} \dot v_j \,.
\label{eq:molecules-from-mixtures}
\ese
Note that if the $i$th row of $S$ is zero, then $x_i$ is constant and we have a redundant equation.
So without loss of generality, we assume that $S$ has no zero rows.

The physical intuition behind a reaction $v_i{\buildrel e_\ell \over \rightarrow }\, v_j$ in a solution
of chemicals is that the reaction rate is proportional to the probability that all the necessary
molecules in $v_i$, the tail of the arrow $e_\ell$, ``meet" in some small volume (this is called the
\emph{mass action principle}). The probability that molecule $r$ is present in some small volume is
proportional to  $x_r$, its concentration in the chemical mix. Assuming these probabilities are
independent of one another, we see that the probability that all the right molecules of $v_i$ are
present in the small volume equals the product of all the concentrations of the molecules in $v_i$.
This product is called $\psi_i(X)$, and these form a vector $\psi(x)$ in $\R^v$.
With the above definition of $S$, we see that this product is proportional to $\prod_j x_j^{S_{ji}}$.
We thus define a vector in vertex space $\R^v$ (using the convention that $0^0:=1$):
\bse
\psi_i(x)\equiv \prod_j x_j^{S_{ji}} \quad \logor \quad \Ln \psi(x)=S^T \Ln x \,.
\label{eq:massaction}
\ese

Next, we transform $\psi\in \R^v$ to the vector in the edge space $\R^{e}$ whose $\ell$th component
is the reaction rate of the $\ell$th reaction $v_i{\buildrel e_\ell \over \rightarrow }\, v_j$.
From the previous paragraph, we conclude that the rate of
the $\ell$th reaction is proportional to $\psi_i(x)$, where the $i$th vertex is the tail (the begin
point) of the $\ell$th directed edge. Thus using the begin matrix $B$ of Section \ref{chap:Laplacians},
we see that the reaction rates are proportional to
\bsenn
B^T\psi(x)\in \R^{e}.
\esenn
In the chemical literature, this proportionality is (nearly) always expressed by a constant
called $k$. It is important to note that this constant is associated with the $\ell$th reaction --
or edge -- and not with the reacting mixture -- or begin vertex -- of that reaction.
Thus we weight the edges using an $e \times e$ diagonal matrix $W$ whose $\ell$th diagonal element
equals a (strictly) positive constant $k_\ell$. The reaction rates are therefore given by
\bsenn
WB^T\psi(x)\in \R^{e}.
\esenn

The reaction $v_i{\buildrel e_\ell \over \rightarrow }\, v_j$ adds to the concentration of mixture
$v_j$ and subtracts from the concentration of mixture $v_i$, both at the rate $k_\ell\psi_\ell$.
Again, with the definitions of $E$ and $B$ of Section \ref{chap:Laplacians}, we compute the rates of
change of the concentration reacting mixtures $v \in \R^v$ as:
\bsenn
\dot v= (E-B)WB^T\psi(x) = \partial WB^T\psi(x) = - L_{\mathrm{out}}^T \psi(x) .
\esenn

Finally, in chemical situations we can't necessarily measure or observe directly the concentrations
of reacting mixtures. Rather, we observe the concentrations of the various molecules $x_i$.
Applying \eqref{eq:molecules-from-mixtures} gives us the final form of the dynamical system in $\R^c$
associated to chemical reaction networks
\bse
\dot x = - S L_{\mathrm{out}}^T \psi(x) \,.
\label{eq:CRN1A}
\ese
Solutions of this system can also be derived from the solutions of the following system:
\bse
\dot v = - L_{\mathrm{out}}^T \psi(Sv)\in \R^{v},
\label{eq:CRN2}
\ese
where we used \eqref{eq:molecules-from-mixtures}.
Interestingly, the reverse is \emph{not} necessarily true. A solution of \eqref{eq:CRN1A} does
\emph{not} always determine a unique solution of \eqref{eq:CRN2}. In fact, one of the problems that
comes up in this type of system, is whether non-trivial reactions can take place even though $\dot x=0$.
From the above equations one can see that could happen if during these reactions $\dot v\in \ke S$.
This is of course impossible if $\ke S\cap\im L_{\mathrm{out}}^T=0$, as we will see in Section
\ref{chap:0defic}.

To summarize the whole framework schematically, here is a diagram of the transformations involved in \eqref{eq:CRN1A}.
{\large \bse
\R^c \;\buildrel S \over \longleftarrow \; \R^v \; \underbrace{\buildrel \partial \over \longleftarrow \; \R^e \;
\buildrel W \over \longleftarrow \; \R^e \; \buildrel B^T \over \longleftarrow}_{-L_{\mathrm{out}}^T}
\; \R^v \; \buildrel \psi \over \longleftarrow \; \R^c.
\label{eq:CRN3}
\ese }

\vskip0.0in
The important step here is that we split these transformations into a non-linear part $\psi(x)$
and a linear part $-S L_{\mathrm{out}}^T$. In the literature, however, since the revolutionary work done in
the 1970's \cite{horn_1972, horn_jackson_1972, feinberg_1972}, the traditional split in treatment
has been between $S\partial$ on the one hand and $W B^T \psi$ on the other. This was done, because
the weights in $W$ are the reaction rates and these are notoriously difficult to measure. In addition,
of course, one did not have access to Theorems \ref{thm:rightkernel} and \ref{thm:leftkernel}.
And so some of the linear transformations --- to wit: $WB^T$ in \eqref{eq:CRN3} --- were lumped
with the non-linear part $\psi$.
What we exhibit in this work is the price paid for that choice.

\vskip0.0in
Our next result is a reality check. Since concentrations cannot be negative, we want to make sure that
the set $\R^c_+=\{x\in\R^c\,|\,\forall i\;:\;x_i\geq 0\}$, also called the positive orthant, is forward invariant.

\vskip0.0in\noindent
\begin{prop} The positive orthant is forward invariant under the flow of \eqref{eq:CRN1A}.
\label{prop:orthant}
\end{prop}

\vskip0.0in\noindent
{\bf Proof.} Suppose there is an orbit $x(t)$
of the flow defined by \eqref{eq:CRN1A} that leaves the positive orthant. Let us say, for some
$\epsilon>0$, $x_j(t_1)=\epsilon$ and $x_j(t_2)=-\epsilon$ crossing the plane $x_j=0$ at the point $P$.
Then by continuity, all orbits with initial condition in some (small) neighborhood $N$ of $P$ of the
plane $x_j=0$, will leave the positive orthant. Thus the \emph{flux} must satisfy
\bsenn
\int_N\,\left(- S L_{\mathrm{out}}^T \psi(x)\cdot \hat e_j\right) \, dA<0 \,,
\esenn
where $(,)$ denotes the standard inner product, $\hat e_j$ is the unit normal to $x_j=0$ pointing
\emph{into} the positive orthant, and $dA$ is the standard ($v-1$)-dimensional area form.
To get the contradiction, it is therefore sufficient to show that if $x_j=0$, then
$\left(- S L_{\mathrm{out}}^T \psi(x)\right)_j\geq 0$.

So suppose $x_j=0$. Since $S$ has no zero rows, there must be a $i$ such that $S_{ji}$ is a positive
integer. From \eqref{eq:massaction} we see that for all $i$ such that then $S_{ji}>0$, we have
$\psi_i=0$. The off-diagonal elements of $-L_{\mathrm{out}}^T$ are non-negative, and so for these
same $i$
\bsenn
\left(-L_{\mathrm{out}}^T \psi\right)_i=\sum_j\,\left(-L_{\mathrm{out}}^T\right)_{ij}\psi_j \geq 0 .
\esenn
Using again that $S_{ji}$ is non-negative, we have
\bsenn
-\left(S L_{\mathrm{out}}^T \psi\right)_j= \sum_i\,S_{ji} \left(-L_{\mathrm{out}}^T
\psi\right)_i\geq 0 .
\esenn
This proves the proposition.
\QED

The preceding development shows that an out-degree Laplacian arises naturally in the analysis of CRNs.
We will see that the algebraic results in Section \ref{chap:Laplacians} are of great use in this
analysis. However to make use of them, we will need to adapt them to the out-degree Laplacian.
Fortunately, this is extremely simple thanks to the dual relationship between the two, namely
$L_{\mathrm{out}}(G)=L_{\mathrm{in}}(G^\Lsh)$ noted in definition \ref{def:laplacians}. Theorems
\ref{thm:multiplicity}, \ref{thm:rightkernel}, and \ref{thm:leftkernel} hold for $L_{\mathrm{out}}$
if we replace each instance of reach and cabal with the dual notions of co-reach and co-cabal.

It turns out that in the development of our theory, we do not use the fact that $S$ is an integer matrix
nor the fact that the Laplacian is out-degree. Hence in the next few sections, we consider the following
slightly more general problem.

\vskip 0.1in\noindent
{\bf Remark:} From now on, the matrix $S$ is a non-negative matrix with no zero rows,
$\psi:\R_+^c\rightarrow \R_+^v$ is defined in \eqref{eq:massaction} and $L$ (the Laplacian) is
$v\times v$ matrix with non-negative diagonal and non-positive off-diagonal elements whose row-sums
all give zero. We consider the system given by \eqref{eq:CRN1A}.

\vskip 0.1in
Equation \eqref{eq:CRN1A} implies that $\dot x \in \im SL_{\mathrm{}}^T$. Thus the orthogonal
projection of $x$ to $\left(\im SL_{\mathrm{}}^T\right)^\perp=\ke L_{\mathrm{}}S^T$ is in
fact a constant of the motion. This motivates the following definition.

\vskip 0.0in\noindent
\begin{defn} Let $P : \R_+^c\rightarrow \ke L_{\mathrm{}}S^T$ be the orthogonal projection.
For $z\in \im P$, let
\bsenn
X_z:=\{x\in \R_+^c\,:\, P(x)=z\} ,,
\esenn
These sets are invariant under the flow of \eqref{eq:CRN1A} and will be referred to as invariant sets.
\label{def:affinespace}
\end{defn}

%\vskip 0.4in
\begin{centering}\section{The Laplacian Zero Deficiency Theorem}
 \label{chap:0defic}\end{centering}
\setcounter{figure}{0} \setcounter{equation}{0}

We present two definition for the deficiency of a network. The fact that they are equal follows from
Proposition \ref{prop:linalg5}.

\begin{defn} The Laplacian deficiency of a chemical reaction network is given by
\bsenn
\delta_L\equiv \dim \left[\ke  S \cap \im  L_{\mathrm{}}^T\right]=
\dim \ke  S  L_{\mathrm{}}^T - \dim \ke  L_{\mathrm{}}^T \,.
\esenn
\label{def:Ldeficiency}
\end{defn}

\vskip-0.3in\noindent
{\bf Remark:} Note that $\delta_L=0$ means that $\ke  S \cap \im  L_{\mathrm{}}^T=\{0\}$, and thus
$\dim \im  S  L_{\mathrm{}}^T = \dim \im  L_{\mathrm{}}^T$.

\vskip0.0in\noindent
This remark and theorems \ref{thm:rightkernel} and \ref{thm:leftkernel} motivate the following convention.

\vskip0.0in\noindent
\begin{defn} Suppose a chemical reaction network (or CRN) has $\delta_L=0$ and its graph has $v$
vertices and $k$ reaches. We will let $\{r_1,\cdots, r_{v-k}\}$ denote a basis of $\im SL^T$.
\label{def:basisImSLT}
\end{defn}

\vskip-0.1in
The next result shows that a 0 deficiency network has a strictly positive equilibrium if and only if
it is CSC. In the two sections that follow we will refine this to show that if a 0 deficiency network
is CSC, then every invariant set $X_z$ (Definition \ref{def:affinespace}) has a unique equilibrium
(Theorem \ref{thm:EandU}) and furthermore, that equilibrium is asymptotically stable (Theorem
\ref{thm:convergence}).

\vskip-0.0in\noindent
\begin{theo}[Laplacian Zero Deficiency Theorem] Suppose a chemical reaction network (or CRN) has $\delta_L=0$.
Then the CRN has a (strictly) positive equilibrium if and only if $G$ is CSC.
\label{thm:0deficiency}
\end{theo}

\vskip-0.1in\noindent
{\bf Proof.} We first prove $\Longrightarrow$. From equation \eqref{eq:CRN1A} we see that
the existence of a positive equilibrium together with $\delta_L=0$ implies that there is a positive
vector $\psi^*=\psi(x^*)$ such that $L_{\mathrm{}}^T\psi^*=0$. From Theorem \ref{thm:leftkernel},
we conclude that (recalling that $k$ is the number of reaches)
\bse
\psi^*= \sum_{m=1}^k a_m \bar\gamma_m^T,
\; \logand \forall \, m,\,a_m>0.
\label{eq:nullvector}
\ese
Furthermore, since $x^*>0$, we have $\psi^*>0$ and so from the form of the $\bar\gamma_m$,
one notes that each reach must be a cabal, and thus (see remarks after Definition \ref{def:CSC}) a strong component. Thus $G$ is CSC.

Now we prove $\Longleftarrow$. Suppose that every reach is a strong component, then using
$\delta_L=0$ we must show that \eqref{eq:nullvector} has a positive solution $\psi^*$ with $x^*>0$.
By positivity, we can take the componentwise logarithm of both sides.
We note that $\mathrm{Ln}\,\psi(x^*)=S^T \mathrm{Ln}\,x^*$. The logarithm of the right hand
side of \eqref{eq:nullvector} can be written as
\bsenn
\mathrm{Ln }\,\sum_{m=1}^k a_m \bar\gamma_m^T = \sum_{m=1}^k (\ln a_m )\, \mathbf{1_{R_m}}+
\mathrm{Ln}\,\sum_{m=1}^k \bar\gamma_m^T ,
\esenn
where $\mathbf{1_{R_m}}$ is the characteristic vector of the $m$th reach or (in this case) component.
Note that
$\sum_{m=1}^k \bar\gamma_m$ has all components positive by assumption. Thus from \eqref{eq:nullvector}
we see that we need to solve $x^*$ in
\bse
S^T \mathrm{Ln}\,x^* = \sum_{m=1}^k (\ln a_m )\, \mathbf{1_{R_m}} +
\mathrm{Ln}\,\sum_{m=1}^k \bar\gamma_m^T .
\label{eq:solvethis}
\ese
This can be re-arranged as
\bse
\mathrm{Ln}\,\sum_{m=1}^k \bar\gamma_m^T =
S^T \mathrm{Ln}\,x^* - \sum_{m=1}^k (\ln a_m )\, \mathbf{1_{R_m}}.
\label{eq:re-arranged}
\ese
We observe that the first term of the right hand side ranges over $\im S^T$ and the second
over $\ke L_{\mathrm{}}$. This has a solution if
\bsenn
\im S^T + \ke L_{\mathrm{}}=\R^{v}.
\esenn
However, this is guaranteed by applying Proposition \ref{prop:linalg4} to the zero deficiency condition.
\QED

\vskip-0.1in
This is the analogue of the classical zero deficiency theorem. It can, however, be strengthened significantly
with very little effort. Here we first show that once can weaken the existence of an positive equilibrium to
the existence of an orbit $x(t)$ such that $\Ln x(t)$ is bounded. Secondly, even if the zero deficiency system
is not CSC, we can still write down equations that determine all the equilibria of the dynamics in each
$X_z$ of Definition \ref{def:affinespace}. In Sections \ref{chap:uniqueness} and \ref{chap:Lyapunov},
we will furthermore show existence and uniqueness as well as asymptotic stability of these equilibria.

\vskip-0.0in\noindent
\begin{theo} Suppose a chemical reaction network (or CRN) has $\delta_L=0$. Then the  CRN has an orbit
$x(t)>0$ such that $\ln x_i(t)$ is bounded for all $i$ if and only if $G$ is CSC.
\label{thm:0deficiencybdd}
\end{theo}

\vskip0.0in\noindent
{\bf Proof.} $\Longleftarrow$ follows from Theorem \ref{thm:0deficiency}. For the other direction, we
compute
\bse
\frac{x(\tau)-x(0)}{\tau} = \frac1\tau \int_0^\tau\dot x\,dt =
-\frac1\tau \int_0^\tau S L_{\mathrm{}}^T \psi(x(t))\,dt=
-S L_{\mathrm{}}^T \frac1\tau \int_0^\tau \psi(x(t))\,dt .
\label{eq:xbounded}
\ese
The requirement on $x_i$ implies that $F(\tau):=\frac1\tau \int_0^\tau \psi(x)\,dt$ has a compact
range of the form $[\epsilon,\epsilon^{-1}]$ for some $\epsilon>0$. Thus $F(n)$ must have a
subsequence $F(n_i)$ convergent to some $F_\infty>0$. On the other hand, the boundedness of $x$
ensures that left hand side of \eqref{eq:xbounded} converges to 0 as $\tau$ tends to infinity. Thus
for the subsequence $\{n_i\}$
\bsenn
0=\lim_{i\rightarrow \infty}\,\frac{x(n_i)-x(0)}{n_i} = -S L_{\mathrm{}}^T F_\infty .
\esenn
The remainder of the proof is as in the first part of Theorem \ref{thm:0deficiency} with $F_\infty$ replacing $\psi^*$.
\QED

Two comments are in order here. The first is that Theorems \ref{thm:0deficiency} and \ref{thm:0deficiencybdd}
imply that for a deficiency zero system $S$ with associated graph $G$ the following holds:
\bsenn
G \textrm{ is CSC} \Longleftrightarrow S \textrm{ has equilibrium}
\Longleftrightarrow S \textrm{ admits orbit $x$ with $\Ln x$ bounded} \,.
\esenn
In particular, for a (Laplacian) deficiency zero system, we have that that CSC implies that none of the
concentrations $x_i$ tend to zero. The persistence conjecture \cite{feinberg_1987} says that this is true
independently of the deficiency.

\vskip-0.0in\noindent
\begin{theo} Suppose a chemical reaction network (or CRN) has $\delta_L=0$ and its underlying graph has
$v$ vertices and $k$ reaches. Then the equilibria in $X_{z_0}$ (see Definition \ref{def:affinespace})
must satisfy these $v$ equations in $v$ unknowns (the $u_i$ and $a_i$):
\bsenn
\psi\left(z_0+\sum_{i=1}^{v-k}\, u_ir_i\right)= \sum_{m=1}^k\, a_m \bar \gamma_m^T \,,
\esenn
where the $r_i$ are and $\bar \gamma_m$ are given in Definition \ref{def:basisImSLT} and Theorem
\ref{thm:leftkernel}.
\label{thm:0deficiencyeqn}
\end{theo}

\vskip-0.1in\noindent
{\bf Proof.} For deficiency zero systems, $x^*$ is an equilibrium if and only if
$\psi(x^*) \in \ke L^T$. Using the basis of Theorem \ref{thm:leftkernel}, this reads
\bsenn
\psi(x^*)=\sum_{m=1}^k\, a_m \bar \gamma_m^T \,.
\esenn
Since $X_{z_0}$ is given by $\left\{z_0+\sum_{i=1}^{v-k}\, u_ir_i\mid z_0\in \ke LS^T, u_i \in \R\right\}$,
the statement follows.
\QED

%\vskip 0.4in
\begin{centering}\section{Existence and Uniqueness of Equilibria}
 \label{chap:uniqueness}\end{centering}
\setcounter{figure}{0} \setcounter{equation}{0}

We will show that a CRN whose associated graph is CSC with zero Laplacian deficiency has \emph{exactly}
one positive equilibrium in each invariant set $X_z$ (see Definition \ref{def:affinespace}).
Informally speaking, then, the set of equilibria forms a graph (as in: \emph{is a function of}) over
$\ke L_{\mathrm{}}S^T$. The precise formulation is given below in Theorem \ref{thm:EandU}.
(See the last paragraph of the introduction for the notation.)

\begin{lem} Given a CSC system with $\delta_L=0$. Suppose $x^*>0$ is an equilibrium.
Then $x>0$ is an equilibrium iff
$\Ln [\psi(x)/ \psi(x^*)]\in \ke L_{\mathrm{}}$, which is equivalent to
$\Ln [x/x^*]\in \ke L_{\mathrm{}}S^T$.
\label{lem:relation-between-eqla}
\end{lem}

\vskip-0.1in\noindent
{\bf Proof.} By our hypotheses, $x^*>0$ is an equilibrium iff $\psi(x^*)= \sum_{i=1}^{k}a_i\bar\gamma_i^T$,
with all $a_i>0$. Similarly, the fact that $x>0$ is an equilibrium is equivalent to
$\psi(x)= \sum_{i=1}^{k}b_i\bar\gamma_i^T$, with all $b_i>0$.
Thus, given that $x^*>0$ is an equilibrium, the same holds for $x$ iff
\begin{align*}
\psi(x)/ \psi(x^*)= \sum_{i=1}^{k}\frac{b_i}{a_i} \mathbf{1_{R_i}} \;\Longleftrightarrow\;
\Ln [\psi(x)/ \psi(x^*)]= \sum_{i=1}^{k}\ln \frac{b_i}{a_i} \mathbf{1_{R_i}} \;\Longleftrightarrow\;
S^T\Ln [x/x^*] = \sum_{i=1}^{k}\ln \frac{b_i}{a_i} \mathbf{1_{R_i}} \,,
\end{align*}
where we used that $\Ln \psi = S^T \Ln x$. Using Theorem \ref{thm:rightkernel}, we get $LS^T\Ln [x/x^*]=0$,
implying the lemma. \QED

\vskip-0.2in\noindent
\begin{prop} Given a CSC system with $\delta_L=0$. For every $z\in \ke L_{\mathrm{}}S^T$, there
exists $y\in \im SL_{\mathrm{}}^T$ such that $x=y+z$ is a positive equilibrium.
\label{prop:existence}
\end{prop}

\vskip-0.0in\noindent
{\bf Proof.} By Theorem \ref{thm:0deficiency}, we may fix a positive equilibrium $x^*$. We also
fix $z\in \ke L_{\mathrm{}}S^T$. By Lemma \ref{lem:relation-between-eqla}, $x$ is a positive equilibrium
if (using the componentwise multiplication $\odot$) it can be written as $x=x^*\odot \Exp \mu$ with $\mu\in \ke L_{\mathrm{}}S^T$. Thus it is
sufficient to show that there is a $\mu^*\in \R^c$ so that $y:=(x^*\odot \Exp \mu^*-z)$ is orthogonal to
$\ke L_{\mathrm{}}S^T$ (i.e. is in $\im SL_{\mathrm{}}^T$), for then $x=z+y$ is a positive equilibrium.
Thus we wish to prove that given $x^*$ and $z$,
\begin{align}
\exists\;\mu^*\in \ke L_{\mathrm{}}S^T \;\; \such \;\; \forall v\in \ke L_{\mathrm{}}S^T\;:\;\;
\left(x^*\odot \Exp \mu^*-z,v\right) = 0 ,
\label{eq:Dg=0}
\end{align}
where $(,)$ stands for the usual inner product. We settle this by defining a smooth function $g(\mu)$ whose
gradient $\nabla g$ with respect to $\mu$ equals $x^*\odot \Exp \mu-z$ and which has a minimum at $\mu^*$,
so that $\nabla g(\\mu^*)=0$.

\begin{figure}
\center
\includegraphics[width=3.5in]{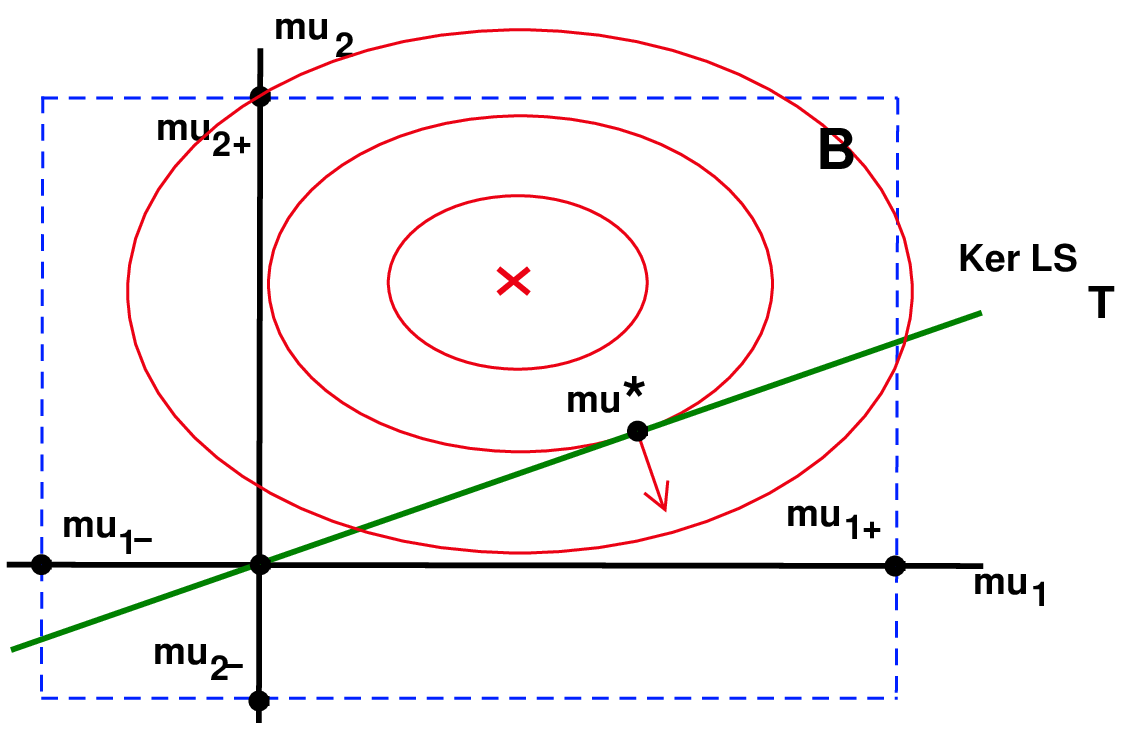}
\caption{\emph{$g(\mu)=(x_1e^{\mu_1}-z_1\mu_1)+(x_1e^{\mu_1}-z_1\mu_1)$. $g(\mu)>\sum_{i=1}^{c} x^*_i$
outside the box and $g(0)=\sum_{i=1}^{c} x^*_i$. Therefore $g$ restricted to $\ke LS^T$ has a
minimum inside the box B.}}
\label{fig:min-gmu}
\end{figure}
To accomplish this, fix $x^*$ and $z$ as above and define $g:\R^{c}\rightarrow \R$ by setting
\begin{align*}
g(\mu)= \left(x^*,\Exp \mu\right)-\left(z,\mu\right) .
\end{align*}
This function is a sum of $c$ one-dimensional functions described in Lemma \ref{lem:calc3}. This Lemma
implies that there is a box $B=[\mu_{1,-},\mu_{1,+}]\times \cdots \times [\mu_{c,-},\mu_{c,+}]\subseteq \R^c$
so that
\begin{align*}
\forall \mu \not\in B\;:\;\; g(\mu) > \sum_{i=1}^{c} x^*_i = g(0) \,.
\end{align*}
See Figure \ref{fig:min-gmu}. Therefore the set $C$ defined by
\begin{align}
C\equiv \left\{\mu\in \ke L_{\mathrm{}}S^T\;\big|\; g(\mu)\leq g(0)\right\}
\end{align}
is non-empty (as it contains 0), closed (by continuity of $g$), and bounded.

Now we restrict $g$ to $\ke L_{\mathrm{}}S^T$. Since $0\in \ke L_{\mathrm{}}S^T$, the
continuous function $g$ assumes its minimum in $\ke L_{\mathrm{}}S^T$ at a point $\mu^*$.
Since $g$ is also differentiable, at $\mu=\mu^*$, we must have
\begin{align*}
\forall v\in \ke L_{\mathrm{}}S^T\;:\;\; 0 = \left(\nabla g(\mu^*),v\right) = \lim_{\epsilon\rightarrow 0}
\dfrac{g(\mu^*+\epsilon v)-g(\mu^*)}{\epsilon}=\left(x^*\odot \Exp \mu^*-z,v\right) ,
\end{align*}
which establishes equation \eqref{eq:Dg=0}, thereby proving the proposition.
\QED

\vskip-0.1in\noindent
\begin{prop} Given a CSC system with $\delta_L=0$.
For every $z\in \ke L_{\mathrm{}}S^T$, there exists \emph{at most one}
$y\in \im SL_{\mathrm{}}^T$ such that $x=y+z$ is a positive equilibrium.
\label{prop:uniqueness}
\end{prop}

\vskip-0.0in\noindent
{\bf Proof.} Suppose that we have $y$ and $u$ both satisfying the requirements. Then by Lemma
\ref{lem:relation-between-eqla},
\begin{align*}
\Ln(z+y)-\Ln(z+u) \in \ke L_{\mathrm{}}S^T ,
\end{align*}
and by hypothesis $y$ and $u$ in $\im SL^T$, so
\begin{align*}
(z+y)- (z+u) \in \im SL_{\mathrm{}}^T .
\end{align*}
By Lemma \ref{lem:linalg2}, the two are orthogonal. Taking the inner product of the two differences gives
\begin{align*}
\left((z+y)- (z+u),\Ln(z+y)-\Ln(z+u)\right)=0 .
\end{align*}
Lemma \ref{lem:calc2} then shows that $z+y=z+u$, and therefore $y=u$.
\QED

\vskip-0.1in\noindent
Putting the last two propositions together immediately gives the main result of this section.

\vskip-0.0in\noindent
\begin{theo} For a CSC system with $\delta_L=0$ we have the following. For every
$z\in \ke L_{\mathrm{}}S^T$, there is a unique
$y\in \left(\ke L_{\mathrm{}}S^T\right)^\perp=\im SL_{\mathrm{}}^T$ such that $x=y+z$ is a positive
equilibrium.
\label{thm:EandU}
\end{theo}

\vskip-0.0in\noindent
{\bf Proof.} Proposition \ref{prop:existence} proves existence and Proposition \ref{prop:uniqueness}
proves uniqueness.
\QED

\vskip-0.2in\noindent
\begin{defn} Suppose $G$ is CSC with $\delta_L=0$.
The unique positive equilibrium of the flow of \eqref{eq:CRN1A} restricted to the invariant sets $X_z$ (see Definition \ref{def:affinespace}) will be denoted by $x^*_z$.
\label{def:unique-equil}
\end{defn}

%\vskip 0.4in
\begin{centering}\section{Convergence to Equilibria}
 \label{chap:Lyapunov}\end{centering}
\setcounter{figure}{0} \setcounter{equation}{0}

For the definition of Lyapunov functions and their use, we refer the reader
to Definitions \ref{def:lyapunov} and \ref{def:omegalimit} and Theorem \ref{thm:lyap}.

The existence of Lyapunov functions depends crucially on the following remarkable result. In the following
proposition and proof, we often refer to $\psi(x^*)$ where $x^*$ is an equilibrium. To avoid cluttering
the formulas, we abbreviate $\psi(x^*)$ as $\psi^*$.

\begin{prop} Let $L$ be the in-degree or out-degree Laplacian. Suppose that there is a $\psi^* > 0$
so that  $\psi^{*T}L=0$. Then the associated graph is CSC and for all $\psi>0$
\bsenn
\psi^T L \left(\Ln \psi-\Ln \psi^*\right)=(\psi-\psi^*)^T L \left(\Ln \psi-\Ln \psi^*\right)\geq 0.
\esenn
Equality holds if and only if on every strong component $C_i$ there is a constant $c_i>0$ such that
\bsenn
\psi|_{C_i}=c_i \psi^*|_{C_i}\,.
\esenn
\label{prop:remarkable}
\end{prop}

\vskip-0.3in\noindent
{\bf Proof.} Let $\psi^* > 0$ and $\psi^{*T}L=0$. We start by observing that Theorem \ref{thm:leftkernel}
implies that then every vertex is in a cabal and so the associated graph is CSC (see Definition
\ref{def:CSC} and the remark following it). Now we write
\bsenn
\psi^T L \left(\Ln \psi-\Ln \psi^*\right)=\psi^T L \,\Ln (\psi/\psi^*)
\esenn
in terms of a sum over its edges. For every directed edge, let $w_{ij}$ be the weight of the edge
$j\rightarrow i$ if $L$ is in-degree Laplacian, and $i\rightarrow j$ if $L$ is an out-degree Laplacian.
Denote by $\sum_{edges}$ the sum over all directed edges. We obtain that $\psi^T L \,\Ln (\psi/\psi^*)$
equals
\benn
\sum_{edges}\,w_{ij}\psi_i\left(\ln \psi_i/\psi^*_i-\ln \psi_j/\psi^*_j\right) &=&
\sum_{edges}\,\psi^*_iw_{ij}\;\;\psi_i/\psi^*_i\left(\ln \psi_i/\psi^*_i-\ln \psi_j/\psi^*_j\right)\\
&\geq & \sum_{edges}\,\psi^*_iw_{ij}\;\;\left( \psi_i/\psi^*_i- \psi_j/\psi^*_j\right) = \psi^{*T}L\;(\psi/\psi^*).
\eenn
The inequality follows from Lemma \ref{lem:calc1} (plus the fact that all non-zero weights are positive).
By assumption, $\psi^{*T}$ is in the left kernel of $L$, and so the last expression gives zero.

Lemma \ref{lem:calc1} also implies the necessary and sufficient condition for equality. To be precise,
that lemma asserts that the condition for equality in the above formula is that on each edge the
value of $(\psi/\psi^*)$ at the head equals its value at the tail. Therefore $(\psi/\psi^*)$ is constant
(and positive) on every strong component.
\QED

\vskip-0.2in\noindent
\begin{prop} Given a Laplacian $L$ with $\delta_L=0$ that has a strictly positive equilibrium $x^*$.
Then the $\omega$-limit set of a positive point $x$ is bounded and is contained in the union of the boundary
of $\R^c_+$ and the set of positive equilibria.
\label{prop:convergence}
\end{prop}

\vskip-0.0in\noindent
{\bf Remark:} Note that $x^*$ in the proposition is an equilibrium.

\vskip-0.0in\noindent
{\bf Proof.}
We will first show that $V:\R^c_+\rightarrow \R$ in \eqref{eq:Lyapunov} defined by
\bse
V(x)=V(x_1,\cdots,x_c):= \sum_{i=1}^{c}\,\int_{x_{i}^*}^{x_{i}}\,\ln(s/x_{i}^*)\,ds
\label{eq:Lyapunov}
\ese
is a Lyapunov function (Definition \ref{def:lyapunov}) for \eqref{eq:CRN1A}. $V$ is clearly
continuously differentiable. In the interest of brevity, we write $\dot V$ for $(\nabla V(x(t)),\dot x(t))$.
To show that $\dot V\leq 0$, we observe
\benn
\dot V &=& \dot x^T \Ln(x/x^*)=-\left(S L_{\mathrm{}}^T\psi(x)\right)^T \mathrm{Ln}(x/x^*)\\
       &=& -\psi^T L_{\mathrm{}}S^T \mathrm{Ln}(x/x^*)= -\psi^T L_{\mathrm{}}\,\mathrm{Ln}(\psi/\psi^*) =
       -\psi^T L_{\mathrm{}}\,\mathrm{Ln}(\psi-\mathrm{Ln}\psi^*) \,.
\eenn
Note from \eqref{eq:massaction} that $x>0$ implies $\psi(x)>0$. It now follows from the first part of
Proposition \ref{prop:remarkable} that $\dot V(x(t))\leq 0$. The second part of Proposition
\ref{prop:remarkable} says that $\dot V(x)=0$ iff $\Ln[\psi/\psi^*]$ is constant on strong components.
Theorem \ref{thm:rightkernel} implies that then $\Ln[\psi/\psi^*]$ is a right null vector of
$L_{\mathrm{}}$. From Lemma \ref{lem:relation-between-eqla} we conclude that then $\dot V(x)=0$ implies
that $x$ is an equilibrium.

Thus by Theorem \ref{thm:lyap},
the $\omega$-limit set of an initial condition may be unbounded, may contain boundary points of the
orthant, and may contain equilibrium points. We rule out the first possibility (unbounded) by showing
that trajectories are bounded. Each integral in the sum of \eqref{eq:Lyapunov} has the form
$I_{x^*_i}(x_i)=\int_{x^*_i}^{x_i}\,\ln s - \ln x^*_i\,ds$, we have
\bsenn
I_{x^*_i}(x_i)=\int_{x^*_i}^{x_i}\,\ln s - \ln x^*_i\,ds= [s\ln s - s - s\ln x^*_i]_{x^*_i}^{x_i}
=x_i(\ln x_i -\ln x^*_i)-(x_i-x^*_i) \geq 0 \,.
\esenn
The final inequality here follows from Lemma \ref{lem:calc1}. Thus each of the integrals in
\eqref{eq:Lyapunov} is non-negative. Furthermore, on the one hand,
$\dot V\leq 0$ and so $V(x(t))\leq V(x(0))$, and on the other,
\bse
I_{x^*_i}(x_i)=\int_{x^*_i}^{x_i}\,\ln s - \ln x^*_i\,ds = x_i(\ln x_i -\ln x^*_i -1)+x^*_i
\label{eq:componentV}
\ese
tends to infinity if $x_i\rightarrow \infty$. This proves that orbits are bounded.
\QED

\vskip 0.0in\noindent
{\bf Remark:} We have required that the system satisfies $\delta_L=0$. This is used to ensure that for any
positive equilibrium $x^*$, we have $\psi(x^*)=\sum_{i=1}^{k}a_i\bar\gamma_i^T$. If we \emph{start}
with the assumption that there is an equilibrium of that form, the hypothesis $\delta_L=0$ is not
necessary.

\vskip-0.0in\noindent
\begin{theo} Suppose $G$ is CSC with $\delta_L=0$. The unique equilibrium $x^*_z$ in $X_z$ (Definition
\ref{def:unique-equil}) is asymptotically stable in $X_z$. The $\omega$-limit set (Definition
\ref{def:omegalimit}) of any positive initial condition either equals that equilibrium or is a bounded
set contained in the boundary of the positive orthant.
\label{thm:convergence}
\end{theo}

\vskip-0.0in\noindent
{\bf Proof.} Existence and uniqueness of $x^*_z$ in $X_z$ follow from Theorem \ref{thm:EandU}. Given
\emph{any} positive equilibrium $x^*$ and consider the function $V$ in \eqref{eq:Lyapunov}. Thus
$V$ is a sum of integrals as exhibited in \eqref{eq:componentV}:
\bsenn
V(x)=\sum_{i=1}^{c}\,I_{x_i^*}(x_i)\,.
\esenn
Now $I_{x_i^*}(x_i)$ has a minimum 0 achieved at $x_i=x_i^*$. On the boundary $\partial X_z$ of $X_z$, at least one of the $x_i$ must be zero. Thus, since $I_{x_i^*}(0)=x_i^*$, $V$ restricted to $\partial X_z$
is greater than or equal to $\min_i x^*_i$. So if we choose an initial condition $x(0)$ such that
$V(x(0))<\min_i x^*_i$, then the fact that $V(x)$ is a Lyapunov function (see the proof of
Proposition \ref{prop:convergence}) prevents $x(t)$ from having limit points in $\partial X_z$.
Thus by Proposition \ref{prop:convergence}, the $\omega$-limit
set is contained in the set of equilibria. But by Theorem \ref{thm:EandU}, the equilibrium is unique.
Thus every orbit with $V(x(0))<\min_i x^*_i$ converges to the unique equilibrium, which therefore is asymptotically stable.  \QED

\vskip 0.0in\noindent
{\bf Remark:} This proof says that the basin (in $X_z$) of attraction of the equilibrium $x^*_z$
contains the set $\{x\in X_z \,:\, V(x)<\min_i x^*_i\}$ where $V$ is given by \eqref{eq:Lyapunov} and
\eqref{eq:componentV}. This gives us a way to get \emph{some} estimate of the basin of attraction
of $x_z$ in $X_z$. Recall that the global attractor conjecture (see Section \ref{chap:intro}) says that in
this case all positive initial conditions converge to the equilibrium.

%\vskip 0.4in
\begin{centering}\section{Examples}
 \label{chap:examples}\end{centering}
\setcounter{figure}{0} \setcounter{equation}{0}

\begin{figure}[!ht]
\center
\includegraphics[width=3.0in]{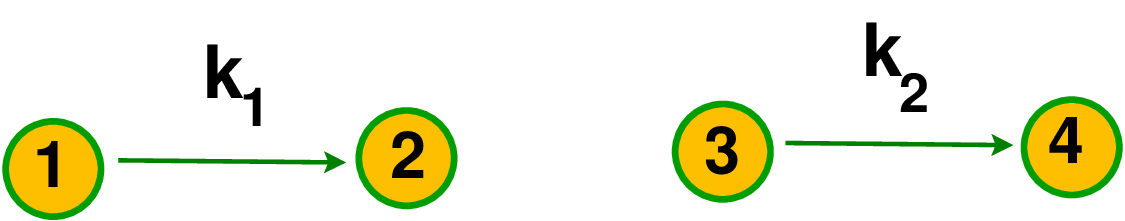}
\caption{\emph{The directed network $G$ of example 1.}}
\label{fig:simple-graph}
\end{figure}

\vskip 0.0in\noindent
{\bf Example 1:} We consider the following simple system taken from wikipedia's ``Chemical reaction network theory" entry:
\benn
\textrm{Reaction 1:}\quad \quad 2\,H_2+ O_2 & \rightarrow & 2\,H_2O \\[-0.0cm]
\textrm{Reaction 2:}\quad\quad \quad C+O_2 &\rightarrow & CO_2
\eenn
The vertices of the network are:
\bsenn
v_1\leftrightarrow 2\,H_2+ O_2 \,,\; v_2\leftrightarrow 2\,H_2O\,,\; v_3\leftrightarrow  C+O_2\,,\;
v_4\leftrightarrow CO_2 \,.
\esenn
The graph $G$ for this system is given in Figure \ref{fig:simple-graph}.
Now we set $x_i$ equal to concentration of following molecules ($[H_2]$ denotes the concentration
 of $H_2$ in chemistry notation).
\bsenn
x_1= [H_2] \,,\; x_2= [O_2]\,,\; x_3= [H_2O]\,,\;
x_4= [C]\,,\;x_5= [CO_2]\,.
\esenn
We assume that all molecules are always mixed uniformly so that the mass action principle applies.
That principle says that in reaction 1, the increase of the number of $x_3$ molecules is proportional
to twice the product of $x_2$ and the square of $x_1$. The (positive) reaction constant is denoted by
$k_1$.  Notice that the increase of $x_3$ molecules must equal the decrease of the $x_1$ molecules.
Reasoning like this we get the following system of equations.
\bse
\begin{matrix}
\dot x_1 & = & -2k_1x_1^2x_2 \\
\dot x_2 & = & -k_1x_1^2x_2 -k_2x_2x_4\\
\dot x_3 & = & 2k_1x_1^2x_2 \\
\dot x_4 & = & -k_2x_2x_4 \\
\dot x_5 & = & k_2x_2x_4
\end{matrix}
\label{eq:simple-example}
\ese

We will now illustrate our methods and main results using this example. Setting up the Laplacian as defined in Section \ref{chap:Laplacians} and $S$ as in Section
\ref{chap:chem}, we get
\bsenn
L_{\mathrm{out}}= \begin{pmatrix} k_1&-k_1&0&0\\0&0&0&0\\0&0&k_2&-k_2 \\0&0&0&0\end{pmatrix} \quad \logand
\quad S = \left( \begin{array}{cccc}
 2&0&0&0\\ 1&0&1&0 \\ 0&2&0&0 \\ 0&0&1&0 \\ 0&0&0&1
 \end{array}\right) \,.
\esenn
One computes
\bsenn
-SL_{\mathrm{out}}^T= \begin{pmatrix} -2k_1&0&0&0\\-k_1&0&-k_2&0\\2k_1&0&0&0\\0&0&-k_2&0\\0&0&k_2&0\end{pmatrix} \quad \logand \quad \Ln \psi= S^T\Ln x= \Ln \begin{pmatrix} x_1^2x_2\\x_3^2\\ x_2x_4 \\x_5\end{pmatrix} \,.
\esenn
Writing out equation \ref{eq:CRN1A}, we obtain \eqref{eq:simple-example} again.

Starting with Section \ref{chap:chem}, one derives with a little computational effort
that the range of $SL_{\mathrm{out}}^T$ is spanned by
\bsenn
\left(\begin{array}{ccccc}1&1/2&-1&0&0\end{array}\right)^T \quad \logand \quad
\left(\begin{array}{ccccc}0& 1 & 0 & 1 & -1\end{array}\right)^T \,,
\esenn
while the kernel of $L_{\mathrm{out}}S^T$ is spanned by
\bsenn
\left(\begin{array}{ccccc}1&0&1&0&0\end{array}\right)^T \quad , \quad
\left(\begin{array}{ccccc}1/2& -1 & 0 & 1 & 0\end{array}\right)^T \quad \logand \quad
\left(\begin{array}{ccccc}-1/2& 1 & 0 & 0 & 1\end{array}\right)^T \,.
\esenn
Definition \ref{def:affinespace} and the remark prior to it now imply that the orthogonal projections
to the latter are preserved by the flow. Thus
\bsenn
c_3=x_1+x_3 \quad , \quad c_4=\tfrac 12 x_1-x_2+x_4 \quad \logand \quad c_5= -\tfrac 12 x_1 +x_2 +x_5
\esenn
are constants of the motion. The dynamics of the system can therefore be described in terms of the variables
\bsenn
u_1=x_1+\tfrac 12 x_2 -x_3 \quad \logand \quad u_2= x_2 +x_4 - x_5
\esenn
plus the constants $c_1$, $c_2$, $c_3$. The conversion is laborious and the resulting equations
are not particularly illuminating, so we leave this as an exercise.

Moving to Section \ref{chap:0defic}, it is not hard to see that $\delta=0$. Since the graph is not
CSC, Theorem \ref{thm:0deficiency} implies that there is no (strictly) positive equilibrium. In this case,
this is reasonably clear from inspecting \eqref{eq:simple-example}. It is even easier to see from the reactions
themselves that eventually some of the substances at the left hand of the reactions must run out.
Since the associated graph (Figure \ref{fig:example-reversed}) is not CSC, Sections \ref{chap:uniqueness}
and \ref{chap:Lyapunov} do not apply.

For the system at hand, we use the above constants of the motion and set $x_1 $ and $x_2$ to be our
independent variables to obtain the equations of Theorem \ref{thm:0deficiencyeqn}. Writing the kernel of
$L_{\mathrm{out}}$ as $\left(\begin{array}{cccc}0& a_1 & 0 & a_2\end{array}\right)^T$ and eliminating
$x_3$, $x_4$, and $x_5$ in favor of the constants $c_i$, we obtain an equilibrium if and only if
\bsenn
\begin{matrix}
x_1^2x_2 & = & 0 \\
(c_3- x_1)^2 & = & a_1\\
x_2(c_4-\tfrac 12 x_1+x_2) & = & 0\\
c_5+\tfrac 12 x_1-x_2 & = & a_2
\end{matrix}
\esenn
And thus given the constants $c_i$, we can solve for $x_1$ , $x_2$, $a_1$, and $a_2$.

\vskip 0.0in\noindent
{\bf Example 2:} Consider the out-degree Laplacian $L_{\mathrm{out}}$ of the graph $G$ in Figure
\ref{fig:example-reversed}. For simplicity, we set all $k_i=1$. The matrix we obtain was given in
equation \eqref{eq:indegree-outdegree}.
This example was chosen to give the same Laplacian as the examples in \cite{kummel, lyons};
its left and right kernels can be found in those papers. In particular, the (right) kernel (Theorem
\ref{thm:rightkernel}) of $L_{\mathrm{out}}$ is spanned by
\begin{align*}
  \gamma_1 = \left(\begin{array}{ccccccc}0& 0 & 1 & 1 & 1& \frac13 &\frac23\end{array}\right)^T \quad
  \logand \quad  \gamma_2={\bf 1}-\gamma_1= \left(\begin{array}{ccccccc}1& 1 & 0 & 0 & 0& \frac23 &\frac13\end{array}\right)^T \,.
\end{align*}
The left kernel (Theorem \ref{thm:leftkernel}) of $L_{\mathrm{out}}$ is spanned by
\begin{align*}
\bar\gamma_1 = \left(\begin{array}{ccccccc}1& 0 & 0 & 0& 0& 0 & 0\end{array}\right) \quad \logand \quad
 \bar\gamma_2= \left(\begin{array}{ccccccc}0& 0&\frac13 & \frac13 & \frac13& 0& 0 \end{array}\right) .
\end{align*}
Let $S$ be given by, for example,
\bsenn
S=\begin{pmatrix} 0&0&3&3&3&1&2\\ 1&1&0&0&0&0&0\\
                  0&2&0&3&0&0&0\\ 0&0&3&0&0&0&2\\
                  0&0&3&0&0&1&0\\ 0&0&0&0&3&0&0\\
\end{pmatrix}
\esenn
One again computes
\bse
-SL_{\mathrm{out}}^T= \begin{pmatrix} 0&0&0&0&0&0&0\\0&0&0&0&0&1&0\\ 0&-2&0&-3&3&0&0 \\ 0&0&-3&3&0&2&-1\\
0&0&-3&3&0&-2&4\\ 0&0&3&0&-3&0&0\end{pmatrix}  \quad \logand \quad
\psi= \begin{pmatrix} x_2 \\ x_2x_3^2 \\ x_1^3x_4^3x_5^3 \\ x_1^3x_3^3 \\ x_1^3x_6^3 \\ x_1x_5 \\ x_1^2x_4^2\end{pmatrix}
\label{eq:psi}
\ese
The evolution equations become
\bse
\begin{matrix} \dot x_1&=&0 \\ \dot x_2&=&x_1x_5\\
\dot x_3 &=& -2x_2x_3^2-3x_1^3x_3^3+3x_1^3x_6^3\\
\dot x_4 &=&-3x_1^3x_4^3x_5^3+3x_1^3x_3^3+2x_1x_5-x_1^2x_4^2 \\
\dot x_5 &=& -3x_1^3x_4^3x_5^3+3x_1^3x_3^3-2x_1x_5+4x_1^2x_4^2 \\
\dot x_6 &=& 3x_1^3x_4^3x_5^3-3x_1^3x_6^3
\end{matrix}
\label{eq:x1constant}
\ese
Note that the second of these equations implies that there is no positive equilibrium at all!
It is possible to show directly that the positive orthant is invariant, but it is much more
involved than in the previous example.
The kernel of the matrix $L_{\mathrm{out}}S^T$ is spanned by $(1,0,0,0,0,0)^T$ and so in this example
the only linear conserved quantity is the value of $x_1$. It is clear from equation
\eqref{eq:x1constant} that it is conserved, though it would take some work to directly verify that
there are no other linear ones.

One confirms (by tedious computation or using symbolic manipulator like MAPLE) that
$\ke  S  L_{\mathrm{out}}^T = \ke  L_{\mathrm{out}}^T$. Thus the Laplacian
deficiency (Definition \ref{def:Ldeficiency}) of this system is zero. There is no strictly positive
equilibrium, and Theorem \ref{thm:0deficiency} says that in this is equivalent
to $G$ not being CSC. This can be directly verified from Figure \ref{fig:example-reversed}.
Better yet, Theorem \ref{thm:0deficiencybdd} implies that there is no orbit such that for all $i$,
$\ln x_i(t)$ is bounded. The orbit of a positive initial condition must approach the boundary of the
orthant or infinity (or both).

In fact, we can use Theorem \ref{thm:0deficiencyeqn} to find the equilibria. Let $\{e_i\}$ denote
the standard basis of $\R^v$. Since $\ke L_{\mathrm{out}}S^T$ is spanned by $e_1$, we can choose
$\{r_2=e_2,\cdots, r_7=e_7\}$ as its orthogonal complement. Setting $x_1=c$ (constant)
and using the above expressions for $\bar\gamma_1$ and $\bar\gamma_2$, the equations for the
equilibria become:
\bsenn
\begin{pmatrix} u_2\\u_2u_3^2\\c^3u_4^3u_5^3\\c^3u_3^3\\c^3u_6^3\\cu_5\\cu_4^2\end{pmatrix}
= \begin{pmatrix} a_1\\0\\a_2/3\\a_2/3\\a_2/3\\0\\0 \end{pmatrix} \,.
\esenn
Let us assume that $c>0$, Then $u_4=u_5=0$. It follows that $a_2=0$, and therefore $u_3=u_6=0$.
The solutions are $u_1=c$ and $u_2=a_1$. One checks directly from \eqref{eq:x1constant} that
$x=(c,d,0,0,0,0)^T$ are indeed equilibria.

One may object that we have overly simplified by setting all the $k_i$ in Figure
\ref{fig:example-reversed} equal to 1. However, the general conclusions are independent of the $k_i$.
The Laplacian for the general case is
\bsenn
L_{\mathrm{out}}= \left[ \begin {array}{ccccccc} 0&0&0&0&0&0&0\\ \noalign{\medskip}-
{\it k_1}&{\it k_1}&0&0&0&0&0\\ \noalign{\medskip}0&0&{\it k_5}&0&-{\it k_5
}&0&0\\ \noalign{\medskip}0&0&-{\it k_3}&{\it k_3}&0&0&0
\\ \noalign{\medskip}0&0&0&-{\it k_4}&{\it k_4}&0&0\\ \noalign{\medskip}
-{\it k_2}&0&0&0&0&{\it k_2}+{\it k_8}&-{\it k_8}\\ \noalign{\medskip}0&0&
-{\it k_6}&0&0&-{\it k_7}&{\it k_7}+{\it k_6}\end {array} \right] \,.
\esenn
Performing the same computations, one shows that the deficiency is still zero, there is one linear constant
of the motion, and the equations for the equilibria can still be written out. The main difference is that
the one linear constant of the motion now cannot easily be read off from the differential equations, because it
depends in a fairly complicated way on the $k_i$. This, in turn, complicates the form of the equations
for the equilibria. Nonetheless, all this can be computed easily using a symbolic manipulator.

Section \ref{chap:uniqueness} and \ref{chap:Lyapunov} assume that $G$ is CSC, and so these have no
further implications for this particular example.

%\vskip 0.4in
\begin{centering}\section{Comparison with Classical Results}
 \label{chap:Comparison1}
 \end{centering}
\setcounter{figure}{0} \setcounter{equation}{0}

We briefly compare our formulation of the main results concerning zero deficiency systems
--- Theorems \ref{thm:0deficiency}, \ref{thm:0deficiencybdd}, \ref{thm:0deficiencyeqn}, \ref{thm:EandU},
and \ref{thm:convergence} --- with their formulation in the literature and show
that some of our results are \emph{strictly} stronger than their classical counterparts.
For this we briefly return to the context of actual chemical reactions (and to the notation
$L_{\mathrm{out}}$ for the Laplacian). Recall the equation \eqref{eq:CRN1A}, governing this type of system
\bsenn
\dot x = - S L_{\mathrm{out}}^T\psi(x)= S\partial\,W B^T \psi(x) .
\esenn
The only nonlinear term is the function $\psi$. So the split in treatment between it and the
linear terms seems very reasonable. However, as explained by equation \eqref{eq:CRN3}, the traditional
split in treatment has been between $S\partial$ on the one hand\footnote{The matrix $S\partial$ is
called the \emph{stoichiometry matrix} in the literature.} and $W B^T \psi$ on the other.
Thus where we find that $\R^c$ is stratified by invariant affine spaces $X_z$ of Definition
\ref{def:affinespace}), the traditional stratification is by the sets where the projection to
$\left( \im S\partial \right)^T=\ke \partial^T S^T$ is constant. Though both are invariant sets, these
sets are \emph{not} the same! We give examples below.

Summarizing, our Theorems \ref{thm:0deficiency} and Theorems \ref{thm:EandU} and \ref{thm:convergence}
imply the classical results. All we need to do is to make the following replacements:
\benn
\delta_L= \dim \ke  S L_{\mathrm{out}}^T\,- \dim \ke  L_{\mathrm{out}}^T \quad
&\textrm{ becomes }\quad  & \delta = \dim \ke  S \partial\,- \dim \ke  \partial\,\\
\ke L_{\mathrm{out}}S^T \quad \logand \quad \im SL_{\mathrm{out}}^T \quad &\textrm{ becomes }\quad  &
\ke \partial^T S^T \quad \logand \quad \im S\partial\\
X_z \quad &\textrm{ becomes }\quad  & \{x_0+V\,:\, V=\im S\partial\} \,.
\eenn

\vskip-0.1in
The following proposition shows that the orthogonal projection onto $\ke L_{\mathrm{out}}S^T$
gives \emph{as many or more} constants of the motion as the projection onto $\ke \partial^T S^T$
(see Definition \ref{def:affinespace}).

\noindent
\begin{prop} $\ke \partial^T S^T\subseteq \ke L_{\mathrm{out}}S^T$.
\label{prop:USimpliesTHEM1}
\end{prop}

\vskip0.0in\noindent
{\bf Proof.}
This becomes clear once we write $L_{\mathrm{out}}$ in full (Definition \ref{def:laplacians}):
\bsenn
\hskip 1.6in \ke \partial^T S^T\subseteq \ke L_{\mathrm{out}}S^T=\ke BW\partial^T S^T\,.  \hskip 2.0in\QED
\esenn

\vskip-0.2in
In example 2 of Section \ref{chap:examples}, there is a linear conserved quantity, namely $x_1=c$.
As mentioned, this is picked up by our method because $(1,0,\cdots ,0)^T$ spans the kernel of
$\ke L_{\mathrm{out}}S^T$. However, the classical theory does not pick up this constant.
Indeed, one checks that the matrix $\partial=E-B$ is given by
\bsenn
\partial =  \left( \begin{array}{cccccccc}
-1&-1& 0& 0& 0& 0& 0& 0\\
 1& 0& 0& 0& 0& 0& 0& 0\\
 0& 0& 0& 0&-1&-1& 0& 1\\
 0& 0& 0& 0& 0& 1&-1& 0\\
 0& 0& 0& 0& 0& 0& 1&-1\\
 0& 1& 1&-1& 0& 0& 0& 0\\
 0& 0&-1& 1& 1& 0& 0& 0
\end {array} \right)
\esenn
Using the same matrix $S$ as before, one obtains that $\ke\partial^T S^T=\{0\}$. Thus the classical
method does not ``see" this constant of the motion.

The next proposition shows that if $\delta=0$, then $\delta_L=0$. However, in the case
of CSC graphs, the two are equivalent. Thus Theorem \ref{thm:0deficiency} is equivalent to
the traditional zero deficiency theorem for these graphs.

\noindent
\begin{prop} $\delta_L\leq \delta$ with equality if $G$ is CSC.
\label{prop:USimpliesTHEM2}
\end{prop}

\vskip0.0in\noindent
{\bf Proof.}
According to Proposition \ref{prop:linalg4}, we have
\bsenn
\left[\ke  S \cap \im \partial\right]^\bot =
\im  S^T+\ke  \partial^T  .
\esenn
This evidently implies that
\bsenn
\delta= \dim \left[\ke  S \cap \im  \partial\,\right]= v -
\dim\left[ \im  S^T+\ke  \partial^T\right] ,
\esenn
where $v$ is the dimension of the vertex space (see Section \ref{chap:chem}). Similarly, we obtain that
\bsenn
\delta_L= v - \dim\left[ \im  S^T+\ke BW\partial^T\right] ,
\esenn
where we have used that $L_{\mathrm{out}}=-BW\partial^T$. The inequality follows from
$\ke \partial^T\subseteq \ke BW\partial^T$.

When $G$ is CSC, Lemma \ref{lem:dimdelta} and Theorem \ref{thm:rightkernel} imply that in addition the
dimensions of $\ke \partial^T$ and $\ke BW\partial^T= \ke L_{\mathrm{out}}$ are equal. Thus the two
sets must be equal.
\QED

\vskip-0.1in\noindent
Similarly, one proves of course that for CSC graphs, the dimensions of
$\ke L_{\mathrm{out}}S^T$ and $\im SL_{\mathrm{out}}^T$ are the same as those of
$\ke \partial^T S^T$ and $\im S\partial$ and so on.

For non-CSC graphs the situation is different as the following simple
example shows. Consider the stargraph with 3 outgoing edges from the central vertex. It is easy to see that
\bsenn
L_{\mathrm{out}}^T=\begin{pmatrix}3&0&0&0\\-1&0&0&0\\-1&0&0&0\\-1&0&0&0\end{pmatrix} \quad
\partial=\begin{pmatrix}-1&-1&-1\\1&0&0\\0&1&0\\0&0&1\end{pmatrix} .
\esenn
The values on the diagonal of the edge weighting matrix $W$ do not matter, so we can take
them to be 1. Suppose further that
\bsenn
S=\begin{pmatrix}2&1&1&1\end{pmatrix} .
\esenn
Since we have
\bsenn
SL_{\mathrm{out}}^T=\begin{pmatrix}3&0&0&0\end{pmatrix} \quad \logand \quad
S\partial = \begin{pmatrix}-1&-2&-2&-2\end{pmatrix} ,
\esenn
it follows that in this example $\delta_L=0$, while $\delta=2$. The obvious generalization
to the stargraph with $k$ outgoing edges will give $\delta_L=0$, while $\delta=k-1$. This shows
that the traditional deficiency can be made arbitrarily larger than the Laplacian one.

\begin{figure}[!ht]
\center
\includegraphics[width=2.5in]{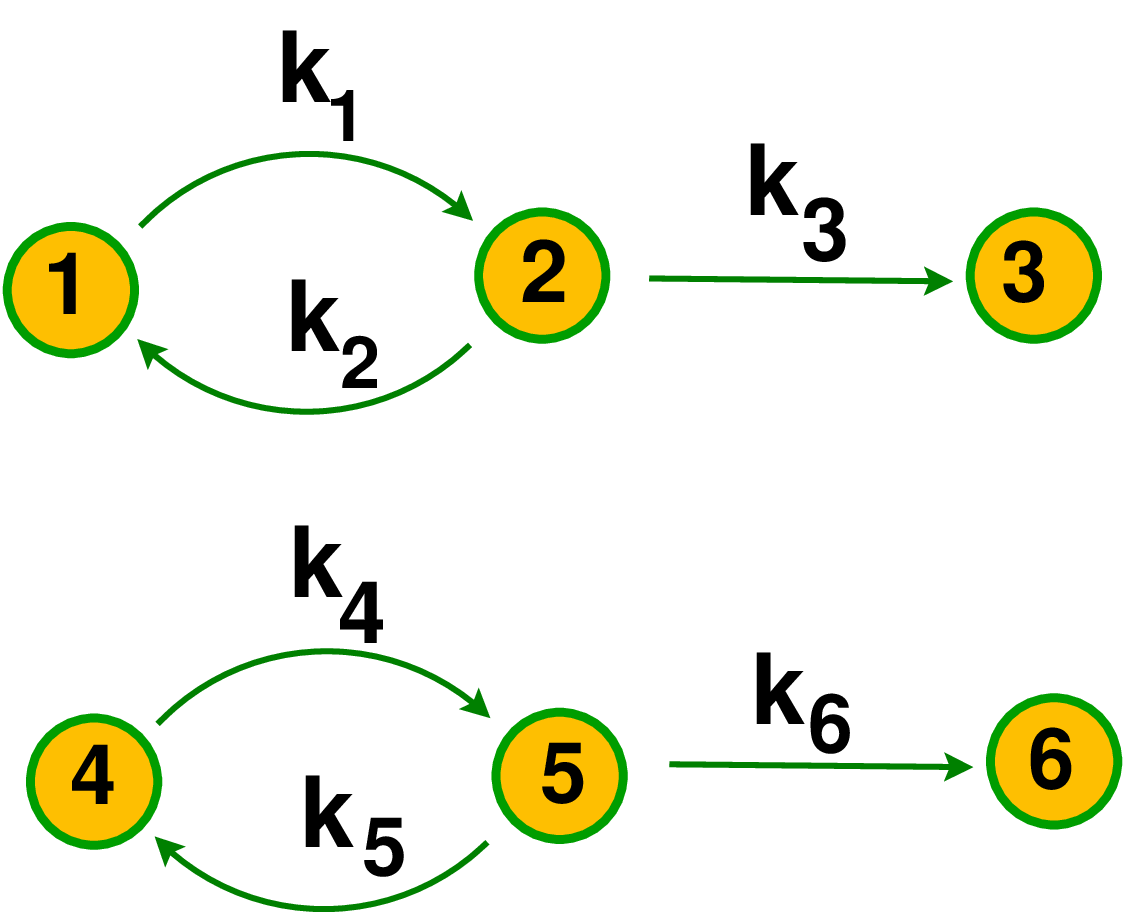}
\caption{\emph{A directed network $G$ with $\delta_L=0$ and $\delta=1$.}}
\label{fig:differentdels}
\end{figure}

\vskip0.0in%\noindent
A much more interesting example of the difference between the Laplacian deficiency $\delta_L$
and the traditional deficiency $\delta$ is an example that plays an important role in the so-called
deficiency one theorem where the additivity of the deficiency is required. We refer to \cite{feinberg_1995}
for the details of that theorem. This example is based on the work \cite{wang_sontag_2008}. Consider the
graph in Figure \ref{fig:differentdels} where the matrix $S$ is given by
\bsenn
S= \left( \begin {array}{cccccc} 1&0&0&0&0&1\\ \noalign{\medskip}1&0&1&0
&0&0\\ \noalign{\medskip}0&1&0&0&0&0\\ \noalign{\medskip}0&0&1&1&0&0
\\ \noalign{\medskip}0&0&0&1&0&1\\ \noalign{\medskip}0&0&0&0&1&0
\end {array} \right) \,.
\esenn
Thus for example, $\psi_3=x_2x_4$ (the 3rd column of $S$) with our usual notation. That means
that the end of the 3rd reaction (or the 3rd vertex) vertex should be labelled $X_2+X_4$, where $X_2$
and $X_4$ represent the molecules for which $x_2$ and $x_4$ represent the concentrations. For convenience,
we list the two relevant matrices.
\bsenn
L_{\mathrm{out}}=\left( \begin {array}{cccccc} {\it k_1}&-{\it k_1}&0&0&0&1
\\ \noalign{\medskip}-{\it k_2}&{\it k_2}+{\it k_3}&-{\it k_3}&0&0&0
\\ \noalign{\medskip}0&0&0&0&0&0\\ \noalign{\medskip}0&0&0&{\it k_4}&-{
\it k_4}&0\\ \noalign{\medskip}0&0&0&-{\it k_5}&{\it k_5}+{\it k_6}&-{\it
k_6}\\ \noalign{\medskip}0&0&0&0&0&0\end {array} \right)
\quad \logand
\esenn
\bsenn
\partial= \left( \begin {array}{cccccc} -1&1&0&0&0&0\\ \noalign{\medskip}1&-1&-
1&0&0&0\\ \noalign{\medskip}0&0&1&0&0&0\\ \noalign{\medskip}0&0&0&-1&1
&0\\ \noalign{\medskip}0&0&0&1&-1&-1\\ \noalign{\medskip}0&0&0&0&0&1
\end {array} \right) \,.
\esenn
If one does the required computations, whose verification we leave as an exercise, it becomes clear
that $\delta_L=0$, while $\delta=1$. Since the graph has two identical weak components, the latter
cannot be additive over weak components. Since that is a condition of the traditional deficiency
one theorem, this system is outside the scope of both theorems. However, it still falls within the scope
of our Laplacian zero deficiency theorem.

\bibliography{ChemReacNets}{}

\begin{thebibliography}{10}

\bibitem{anderson_2011}
D.~F. Anderson.
\newblock A proof of the global attractor conjecture in the single linkage
  class case.
\newblock {\em SIAM Journal on Applied Mathematics}, 71(4):1487--1508, 2011.

\bibitem{aris_1965}
R.~Aris.
\newblock Prolegomena to the rational analysis of systems of chemical
  reactions.
\newblock {\em Archive for Rational Mechanics and Analysis}, 19(2):81--99,
  1965.

\bibitem{bollobas}
B.~Bollob\'as.
\newblock {\em Modern Graph Theory}.
\newblock Springer, 1998.

\bibitem{bray_1921}
W.~C. Bray.
\newblock A periodic reaction in homogeneous solution and its relation to
  catalysis.
\newblock {\em Journal of the American Chemical Society}, 43(6):1262--1267,
  1921.

\bibitem{brunner_2018}
J.~D. Brunner and G.~Craciun.
\newblock Robust persistence and permanence of polynomial and power law
  dynamical systems.
\newblock {\em SIAM Journal on Applied Mathematics}, 78(2):801--825, 2018.

\bibitem{caugh}
J.~S. Caughman and J.~J.~P. Veerman.
\newblock Kernels of directed graph {L}aplacians.
\newblock {\em The Electronic Journal of Combinatorics}, 13(1), 2006.

\bibitem{chebotarev2002}
P.~Chebotarev and R.~Agaev.
\newblock Forest matrices around the {L}aplacian matrix.
\newblock {\em Linear Algebra and its Applications}, 356(1-3):254--273, 2002.

\bibitem{craciun_dickenstein_shiu_sturmfels_2009}
G.~Craciun, A.~Dickenstein, A.~Shiu, and B.~Sturmfels.
\newblock Toric dynamical systems.
\newblock {\em Journal of Symbolic Computation}, 44(11):1551--1565, 2009.

\bibitem{dickenstein_2021}
A.~Dickenstein.
\newblock Algebraic geometry tools in systems biology.
\newblock {\em Notices of the American Mathematical Society}, 67:1706--1715,
  2021.

\bibitem{feinberg_1972}
M.~Feinberg.
\newblock Complex balancing in general kinetic systems.
\newblock {\em Archive for Rational Mechanics and Analysis}, 49(3):187--194,
  1972.

\bibitem{feinberg_1987}
M.~Feinberg.
\newblock Chemical reaction network structure and the stability of complex
  isothermal reactors---i. the deficiency zero and deficiency one theorems.
\newblock {\em Chemical Engineering Science}, 42(10):2229--2268, 1987.

\bibitem{feinberg_1995}
M.~Feinberg.
\newblock The existence and uniqueness of steady states for a class of chemical
  reaction networks.
\newblock {\em Archive for Rational Mechanics and Analysis}, (132):311--370,
  1995.

\bibitem{feinberg_2019}
M.~Feinberg.
\newblock {\em Foundations of chemical reaction network theory}.
\newblock Springer, 2019.

\bibitem{feinberg_horn_1974}
M.~Feinberg and F.~J.~M. Horn.
\newblock Dynamics of open chemical systems and the algebraic structure of the
  underlying reaction network.
\newblock {\em Chemical Engineering Science}, 29(3):775--787, 1974.

\bibitem{gunawardena_2003}
J.~Gunawardena.
\newblock Chemical reaction network theory for in-silico biologists.
\newblock Available online: \url{http://vcp.med.harvard.edu/papers/crnt.pdf},
  2003.

\bibitem{horn_1972}
F.~J.~M. Horn.
\newblock Necessary and sufficient conditions for complex balancing in chemical
  kinetics.
\newblock {\em Archive for Rational Mechanics and Analysis}, 49(3):172--186,
  1972.

\bibitem{horn_jackson_1972}
F.~J.~M. Horn and R.~Jackson.
\newblock General mass action kinetics.
\newblock {\em Archive for Rational Mechanics and Analysis}, 47(2):81--116,
  1972.

\bibitem{joshi_2017}
B.~Joshi and A.~Shiu.
\newblock Which small reaction networks are multistationary?
\newblock {\em SIAM Journal on Applied Dynamical Systems}, 16(2):802--833,
  2017.

\bibitem{kaufman_soule_2019}
M.~Kaufman and C.~Soul{\'e}.
\newblock On the multistationarity of chemical reaction networks.
\newblock {\em Journal of Theoretical Biology}, 465:126--133, 2019.

\bibitem{kim_2018}
Y.~Kim, J.~W. Kim, Z.~Kim, and W.~Y. Kim.
\newblock Efficient prediction of reaction paths through molecular graph and
  reaction network analysis.
\newblock {\em Chemical Science}, 9(4):825--835, 2018.

\bibitem{mirzaev_2013}
I.~Mirzaev and J.~Gunawardena.
\newblock {L}aplacian dynamics on general graphs.
\newblock {\em Bulletin of Mathematical Biology}, 75(11):2118--2149, sep 2013.

\bibitem{pucci_2018}
F.~Pucci and M.~Rooman.
\newblock Deciphering noise amplification and reduction in open chemical
  reaction networks.
\newblock {\em Journal of The Royal Society Interface}, 15(149):20180805, dec
  2018.

\bibitem{schaft_2013}
S.~Rao, A.~J. Van~Der Schaft, and B.~Jayawardhana.
\newblock A graph-theoretical approach for the analysis and model reduction of
  complex-balanced chemical reaction networks.
\newblock {\em Journal of Mathematical Chemistry}, 51(9):2401--2422, 2013.

\bibitem{schaft_rao_jayaw_2015}
A.~J. Van~Der Schaft, S.~Rao, and B.~Jayawardhana.
\newblock A network dynamics approach to chemical reaction networks.
\newblock {\em International Journal of Control}, 89(4):731--745, 2015.

\bibitem{shiu_2010}
A.~Shiu and B.~Sturmfels.
\newblock Siphons in chemical reaction networks.
\newblock {\em Bulletin of Mathematical Biology}, 72(6):1448--1463, jan 2010.

\bibitem{sontag-2001}
E.~D. Sontag.
\newblock Structure and stability of certain chemical networks and applications
  to the kinetic proofreading model of {T}-cell receptor signal transduction.
\newblock {\em Automatic Control, IEEE Transactions on}, 46:1028 -- 1047, 08
  2001.

\bibitem{sontag-2002}
E.~D. Sontag.
\newblock Correction to ``structure and stability of certain chemical networks
  and applications to the kinetic proofreading model of {T}-cell receptor
  signal transduction".
\newblock {\em IEEE Trans Aut Control}, 47(4):1028--1047, 2002.

\bibitem{stern}
S.~Sternberg.
\newblock {\em Dynamical Systems}.
\newblock Dover, 2010, revised 2013.

\bibitem{DG4}
J.~J.~P. Veerman.
\newblock Digraphs {IV}.
\newblock Available online:\\
  \url{http://web.pdx.edu/~veerman/2019-Digraphs-4.pdf}, 2019.

\bibitem{kummel}
J.~J.~P. Veerman and E.~Kummel.
\newblock Diffusion and consensus on weakly connected directed graphs.
\newblock {\em Linear Algebra and its Applications}, 578:184--206, 2019.

\bibitem{lyons}
J.~J.~P. Veerman and R.~Lyons.
\newblock A primer on {L}aplacian dynamics in directed graphs.
\newblock {\em Nonlinear Phenomena in Complex Systems}, 23(2), 2020.

\bibitem{wang_sontag_2008}
L.~Wang and E.~D. Sontag.
\newblock On the number of steady states in a multiple futile cycle.
\newblock {\em J. Math. Biol.}, 1(57):29--52, 2008.

\bibitem{wegsch_1901}
R.~Wegscheider.
\newblock Ueber simultane {G}leichgewichte und die {B}eziehungen zwischen
  {T}hermodynamik und {R}eactionskinetik homogener {S}ysteme.
\newblock {\em Monatshefte fuer Chemie}, 8(32):849--906, 1901.

\bibitem{winfree_1984}
A.~T. Winfree.
\newblock The prehistory of the {B}elousov-{Z}habotinsky oscillator.
\newblock {\em J. Chem. Educ.}, (61):661--663, 1984.

\bibitem{zhabot_1991}
A.~M. Zhabotinsky.
\newblock A history of chemical oscillations and waves.
\newblock {\em Chaos: An Interdisciplinary Journal of Nonlinear Science},
  1(4):379--386, dec 1991.

\end{thebibliography}
\bibliographystyle{plain}

\vspace{\fill}
\end{document}